\documentclass[11pt]{amsart}
\usepackage[T1]{fontenc}
\usepackage{lmodern}
\usepackage[letterpaper,margin=1.15in]{geometry}
\usepackage{amsmath,amssymb,amsthm,mathtools}
\usepackage{booktabs,array}
\usepackage{enumitem}
\usepackage{tikz}
\usetikzlibrary{arrows.meta,positioning,calc,decorations.pathreplacing,patterns}
\usepackage{microtype}
\usepackage[colorlinks=true,linkcolor=blue!45!black,citecolor=blue!45!black,urlcolor=blue!45!black]{hyperref}
\usepackage{xcolor}
\numberwithin{equation}{section}
\newtheorem{theorem}{Theorem}[section]
\newtheorem{proposition}[theorem]{Proposition}
\newtheorem{lemma}[theorem]{Lemma}
\newtheorem{corollary}[theorem]{Corollary}
\theoremstyle{definition}
\newtheorem{definition}[theorem]{Definition}
\newtheorem{question}[theorem]{Question}
\theoremstyle{remark}

\newcommand{\HH}{\mathbb{H}}
\newcommand{\RR}{\mathbb{R}}
\newcommand{\lk}{\operatorname{lk}}
\newcommand{\Aut}{\operatorname{Aut}}
\newcommand{\Vol}{\operatorname{Vol}}
\newcommand{\relint}{\operatorname{relint}}
\newcommand{\inter}{\operatorname{int}}
\newcommand{\arcosh}{\operatorname{arcosh}}
\newcommand{\cA}{\mathcal{A}}
\newcommand{\Tt}{T_2}
\newcommand{\Th}{T_3}
\newcommand{\dd}{\delta}
\newcommand{\ceil}[1]{\left\lceil #1\right\rceil}

\setlist[enumerate]{label=\textup{(\roman*)},leftmargin=*,itemsep=2pt,topsep=4pt}

\begin{document}

\title[On the minimality of the compact right-angled $120$-cell]{On the minimality of the compact right-angled $120$-cell}

\author{Lizi Guo}
\address{Department of Mathematics, University of Wisconsin--Madison, 480 Lincoln Drive, Madison, WI 53706, USA}
\email{llguo@wisc.edu}

\author{Fangting Zheng}
\address{Department of Mathematical Sciences, Xi'an Jiaotong-Liverpool University, Suzhou 215123, China}
\email{fangting.zheng@xjtlu.edu.cn}
\date{October 18, 2026}
\thanks{Fangting Zheng was supported by NSFC No. 12471067.}

\begin{abstract}
We prove that a compact right-angled polytope in hyperbolic $4$-space has at least $120$ facets, with equality only for the right-angled $120$-cell. This answers a question of Potyagailo and Vinberg. As a consequence, the right-angled $120$-cell is the unique compact right-angled $4$-polytope of minimal volume and $f$-vector.
\end{abstract}

\maketitle
\section{Introduction}\label{sec:intro}

A compact polytope $P$ in hyperbolic space $\HH^n$ is \emph{right-angled} if all of its dihedral angles are equal to $\pi/2$. The reflections in the facets of such a polytope generate a right-angled Coxeter group with fundamental domain $P$, and its torsion-free subgroups of finite index give closed hyperbolic $n$-manifolds tessellated by copies of $P$. Right-angled polytopes are therefore a convenient source of explicit hyperbolic manifolds, while at the same time their combinatorics is severely restricted.

In dimension three, by theorems of Pogorelov~\cite{Pog67} and Andreev~\cite{And70} (see also~\cite{HR93,RHD07}), a simple $3$-polytope with more than four faces has a compact right-angled realization in $\HH^3$ if and only if it has no prismatic $3$- or $4$-circuits, and this realization is unique up to isometry. The smallest example is the right-angled dodecahedron $D$, next to which come the L\"obell polyhedra~\cite{Loe31,Ves87}. Inoue described in~\cite{Ino08,Ino22} the first $825$ right-angled polyhedra in order of volume. For surveys of right-angled polyhedra in $\HH^3$ see~\cite{Ves17,Ero19}. Notably, compact right-angled polytopes do not exist in $\HH^n$ for $n\ge5$~\cite{Vin85,PV05}, so $n=4$ is the largest possible dimension. The right-angled $120$-cell $Z$ of volume $34\pi^2/3$~\cite{KP11,MZ21} is the most symmetric example in $\HH^4$. It is combinatorially the regular $120$-cell with $f$-vector $(f_0(Z),f_1(Z),f_2(Z),f_3(Z))=(600,1200,720,120)$, and its facets are right-angled dodecahedra $D$.

Constructing hyperbolic manifolds from Coxeter polytopes has a long history~\cite{Dav85,RT01,CM05,KS16,MZ23}. For the right-angled $120$-cell $Z$, Garrison and Scott~\cite{GS03} investigated some simple small covers of the right-angled dodecahedron $D$ and the right-angled 120-cell $Z$. Ma and Zheng~\cite{MZ21} classified all orientable small covers over the right-angled 120-cell $Z$. As the right-angled 120-cell $Z$ has the smallest volume among all known compact right-angled 4-polytopes, it is natural to ask whether the minimality holds in general. Potyagailo and Vinberg raised this question explicitly~\cite{PV05}:
\begin{quote}
\emph{Is it true that the least number of hyperfaces of a compact right-angled polyhedron in $\HH^4$ is $120$?}
\end{quote}

\noindent We answer this question affirmatively as follows.
\begin{theorem}\label{thm:main}
Every compact right-angled polytope $P\subset\HH^4$ has at least $120$ facets. Equality holds if and only if $P$ is isometric to the right-angled $120$-cell $Z$.
\end{theorem}

By Euler's formula and the Gauss--Bonnet theorem, a compact right-angled polyhedron has at least twelve faces, and the right-angled dodecahedron $D$ is the only one with twelve. Besides, the right-angled dodecahedron $D$ and the L\"obell polyhedron $L_6$ are the two compact right-angled polyhedra of smallest volume~\cite{Ino08}. In dimension four, Kolpakov~\cite{Kol12} proved that the ideal right-angled $24$-cell has the fewest facets and the smallest volume among ideal right-angled polytopes in $\HH^4$. Although the ideal method does not apply to compact right-angled polytopes, Kolpakov observed in~\cite{Kol12thesis} that any compact right-angled polytope $P\subset\HH^4$ with $f_3(P)>120$ has larger volume than the right-angled 120-cell $Z$. He also asked whether a compact right-angled hyperbolic $4$-polytope with a facet other than the right-angled dodecahedron $D$ must have more than $120$ facets and more than $600$ vertices. For a compact right-angled hyperbolic $4$-polytope, its $f$-vector and volume are determined by $f_3(P)$ and the total excess of its $2$-faces over pentagons (Section~\ref{sec:consequences}). Theorem~\ref{thm:main} therefore implies the following, which answers the question of Kolpakov. 

\begin{corollary}\label{cor:volume}
Every compact right-angled polytope $P\subset\HH^4$ satisfies
\[
 f_0(P)\ge600,\qquad f_1(P)\ge1200,\qquad f_2(P)\ge720,\qquad \Vol(P)\ge\frac{34\pi^2}{3}.
\]
Equality in any one of these inequalities holds if and only if $P$ is isometric to $Z$.
\end{corollary}

\begin{corollary}\label{cor:kolpakov}
If a compact right-angled polytope $P\subset\HH^4$ has a facet that is not isometric to $D$, then $f_3(P)\ge121$, $f_0(P)\ge606$ and $\Vol(P)\ge23\pi^2/2$.
\end{corollary}

\subsection*{Outline of the proof}
Let $E$ be the double of the right-angled dodecahedron $D$ across a face, and let $\Tt$ and $\Th$ be the polyhedra obtained by gluing three copies of $D$ in two different ways (Section~\ref{subsec:unions}). Set $\cA=\{D,E,\Tt,\Th\}$. The proof consists of two phases. Let $P$ be a compact right-angled hyperbolic $4$-polytope. First, if $f_3(P)\le120$, then every facet of $P$ belongs to $\cA$. Next, a polytope $P$ all of whose facets belong to $\cA$ can be cut by hyperplanes into copies of the right-angled $120$-cell $Z$.

Let $N=f_3(P)$. The dual boundary complex $K(P)$ of $P$ is a simplicial $3$-sphere which is flag and has no squares. For a facet $A$ of $P$, counting the vertices of $K(P)$ at distance one and two from the vertex $v_A$ dual to facet $A$ gives an identity (Proposition~\ref{prop:identity})
\begin{equation}\label{eq:intro-identity}
 2N=30+5\,m(A)+\sum_{B\sim A}\dd(B,A\cap B)+2\,n_{\ge3}(A,P).
\end{equation}
Here $m(A)$ is the number of faces of polyhedron $A$, $n_{\ge3}(A,P)\ge0$ counts the facets at distance at least three from the polyhedron $A$ in the dual graph, and $\dd(B,A\cap B)$ is a nonnegative integer, the \emph{defect} of a neighbour $B$ of $A$. Only neighbours that are dodecahedra have zero defect. So roughly speaking, $\delta(B,A\cap B)$ measures how far the polyhedron $B$ is from being a right-angled dodecahedron $D$. All the facets of the right-angled 120-cell $Z$ are right-angled dodecahedra $D$, so all defects vanish in the right-angled 120-cell $Z$. By~\eqref{eq:intro-identity}, if $N\le120$, then every facet $A$ satisfies $\dd_P(A)=\sum_{B}\dd(B,A\cap B)\le210-5m(A)$. Intuitively, most neighbours of a facet must be very close to dodecahedra.

Consider polyhedra with fewer than twenty faces. There are only $45$ combinatorial types, and therefore isometry classes, in total, and among them the right-angled dodecahedron $D$ and its double $E$ are the only ones with a regular pentagonal face (Lemma~\ref{lem:small-facts}). Hence a facet $A$ other than $D$ and $E$ has no dodecahedral neighbours, and every neighbour $B$ contributes a positive defect $\dd(B,A\cap B)$. For each pair $(A,B)$ of combinatorial types, we check all possibilities by computer. It turns out that neighbours of small defect are rare and, moreover, that several of them cannot occur simultaneously. Finally, we show in Section~\ref{sec:small} that, when $N\le120$, the only facets with fewer than twenty faces are the right-angled dodecahedron $D$ and its double $E$.

The only neighbours of defect zero are dodecahedra, which can only be attached along regular pentagons. A connected cluster of regular pentagons always lies on the boundary of a single dodecahedron (Lemma~\ref{lem:clusters}). Thus, the number of regular pentagons of a facet polyhedron $A$ is bounded by the number of faces of $A$. Identity~\eqref{eq:intro-identity} then gives $m(A)\le27$ for every facet of $P$ (Proposition~\ref{prop:upper}). A similar estimate shows that a facet with $20\le m(A)\le27$ faces must have many regular pentagons. After a finite search, only $\Tt$ and $\Th$ survive (Section~\ref{sec:classification}). Hence every facet of $P$ belongs to $\cA$ (Proposition~\ref{prop:classification}).

Lastly, as each member of $\cA$ is a union of dodecahedra glued along \emph{internal pentagons}, every side of an internal pentagon lies in a hexagonal or heptagonal face, and we call it a \emph{seam} (Figure~\ref{fig:seams}). When two facets share such a face, they have the same seams on it. For a $4$-polytope $P$, all of whose facets belong to $\cA$, we show that, starting from an internal pentagon and extending it across its seams into the neighbouring facets, we obtain an embedded dodecahedron that is orthogonal to all facets of $P$ it meets. Cutting along such dodecahedra inductively, we dissect $P$ into $q\ge1$ copies of the right-angled $120$-cell $Z$ and obtain $N=106q+14$ (Theorem~\ref{thm:dissection}). If $N\le120$, this forces $q=1$, and $P=Z$.

\subsection*{The role of the computer and data availability}
The computer is used for three types of finite calculation.
\begin{enumerate}
\item The combinatorial types of right-angled polyhedra with at most $27$ faces correspond to the $5$-connected triangulations of the $2$-sphere with at most $27$ vertices. These are generated with \texttt{plantri}~\cite{BM07}. We recheck the lists of $407906$ pairwise nonisomorphic triangulations.
\item For $195$ specific polyhedra, we compute unit normals and edge lengths of their right-angled realizations in ball arithmetic, with a rigorous error bound obtained from a contraction-mapping argument.
\item We enumerate the candidate sets of regular pentagons in Section~\ref{sec:classification}, and the possible ways of attaching a small polyhedron along a common face in Section~\ref{sec:small}.
\end{enumerate}
Appendix~\ref{app:computation} explains computation methods, and Appendix~\ref{app:data} records the data used in the proofs. Programs and data are available at \url{https://github.com/Lizi-Guo/120cell-minimality}, including the candidate enumeration of Section~\ref{sec:classification} and all counts quoted in Sections~\ref{sec:small}--\ref{sec:classification}.

\subsection*{Organization}
Section~\ref{sec:prelim} collects background on right-angled polytopes, their dual complexes and introduces the set $\cA$. Section~\ref{sec:identity} proves the identity~\eqref{eq:intro-identity}. Sections~\ref{sec:small}, \ref{sec:bound} and~\ref{sec:classification} show that a polytope with at most $120$ facets has all its facets in $\cA$. Section~\ref{sec:dissection} carries out the dissection along internal pentagons and completes the proof of Theorem~\ref{thm:main}, and Section~\ref{sec:consequences} deduces Corollaries~\ref{cor:volume} and~\ref{cor:kolpakov} and ends with some questions. The programs are described in Appendix~\ref{app:computation}, and the data they produce are recorded in Appendix~\ref{app:data}.

\subsection*{AI disclosure}
The authors used ChatGPT (OpenAI) to assist in developing ideas, exploring proof strategies, and writing the computational code used in this work. Claude (Anthropic) assisted with drafting the manuscript. The authors take full responsibility for the accuracy and integrity of the paper.

\section{Preliminaries}\label{sec:prelim}

\subsection{Right-angled polytopes}\label{subsec:ra}
We use the hyperboloid model $\HH^n=\{x\in\RR^{n,1}:\langle x,x\rangle=-1,\ x_{n+1}>0\}$, where
\[
 \langle x,y\rangle=x_1y_1+\dots+x_ny_n-x_{n+1}y_{n+1}.
\]
Hyperplanes in $\HH^n$ are the sets $e^\perp\cap\HH^n$ with $e$ a unit spacelike vector. If $F$ is a facet of a polytope $X\subset\HH^n$, we write $H_F$ for the hyperplane spanned by $F$ and $e_F$ for its outward unit normal. Thus $X=\bigcap_F\{x\in\HH^n:\langle x,e_F\rangle\le0\}$, the intersection over all facets $F$ of $X$.
See, for example,~\cite{Rat19} for more details.

We shall use the following facts about a compact right-angled polytope $X\subset\HH^n$ repeatedly. See~\cite[Section~3]{Vin85} and~\cite[Proposition~4.1]{PV05}.
\begin{enumerate}[label=\textup{(F\arabic*)}]
\item If the outward normals of a set $S$ of facets of $X$ are pairwise orthogonal (more generally, if their Gram matrix is positive definite), then the facets in $S$ have a common face of codimension $|S|$.
\item Two facets $F,G$ of $X$ are either adjacent with $\langle e_F,e_G\rangle=0$, or disjoint with $\langle e_F,e_G\rangle<-1$. In the latter case the hyperplanes $H_F,H_G$ are \emph{ultraparallel} and $-\langle e_F,e_G\rangle$ is the hyperbolic cosine of their distance.
\item $X$ is simple, and every face of $X$ of positive dimension is a compact right-angled polytope in its own span.
\end{enumerate}
In particular every $2$-face of $X$ is a compact right-angled polygon with at least five sides, by the Gauss--Bonnet formula.

From now on $P$ always denotes a compact right-angled polytope in $\HH^4$ and $N=f_3(P)$ its number of facets, where $f_i$ is the number of $i$-dimensional faces.
Two facets are \emph{neighbours}, written $A\sim B$, if they share a $2$-face. For a $3$-polyhedron $B$ we write $m(B)=f_2(B)$ for its number of faces, and for a polygon $Q$ we write $f_1(Q)$ for its number of sides.

\subsection{Dual complexes}\label{subsec:dual}
For a simple polytope $X$ let $K(X)$ be the boundary complex of the dual polytope. It is a simplicial $(n-1)$-sphere whose vertices correspond to the facets of $X$, and a set of vertices spans a simplex if and only if the corresponding facets have nonempty common intersection. We write $\Gamma(X)$ for the $1$-skeleton of $K(X)$, the \emph{dual graph} of the facets of $X$. The \emph{distance} between two facets of $X$ is their distance in the graph $\Gamma(X)$, that is, the number of edges of a shortest path joining them. The \emph{degree} $k(F)$ of a facet $F$ of $X$ is the degree of the dual vertex $v_F$ in $\Gamma(X)$. This combinatorial distance should not be confused with the hyperbolic distance between hyperplanes. If $A$ is a facet of $P$ and $v_A$ the corresponding vertex of $K(P)$, then
\begin{equation}\label{eq:link}
 \lk_{K(P)}(v_A)\cong K(A).
\end{equation}

A simplicial complex is \emph{flag} if every set of pairwise adjacent vertices spans a simplex. We say that it has \emph{no squares} if its $1$-skeleton has no induced cycle of length four.

\begin{lemma}\label{lem:FNS}
Let $X$ be a compact right-angled polytope in $\HH^n$. Then $K(X)$ is a flag simplicial $(n-1)$-sphere without squares. In particular, the complex $K(P)$ and all its vertex links have these properties.\qed
\end{lemma}

Here $K(X)$ is exactly the nerve of the Coxeter group corresponding to $X$ when $X$ is right-angled. The following is the uniqueness part of Andreev's theorem~\cite{And70}, in the form in which we shall use it.

\begin{theorem}[Andreev]\label{thm:andreev}
Let $B,B'\subset\HH^3$ be compact right-angled polyhedra. Every isomorphism $\Gamma(B)\to\Gamma(B')$ is induced by an isometry $B\to B'$. In particular, every automorphism of $\Gamma(B)$ is induced by an isometry of $B$, and so preserves edge lengths, the metrics of the faces and the set of regular pentagonal faces.
\end{theorem}

\begin{proof}
By Lemma~\ref{lem:FNS}, $K(B)$ is the clique complex of $\Gamma(B)$, so a graph isomorphism induces an isomorphism of face lattices. Since $B$ has more than four faces (in fact at least twelve, by~\eqref{eq:euler}), the uniqueness part of Andreev's theorem~\cite{And70}, \cite[Theorem~1.4]{RHD07} applies to the two realizations of the same labelled polyhedron with all dihedral angles $\pi/2$, and gives the isometry.
\end{proof}

The combinatorial types we need are therefore encoded by dual graphs. Recall that a triangulation of the $2$-sphere is \emph{$5$-connected} if it has at least six vertices and remains connected after deleting any four of them.

\begin{lemma}\label{lem:5conn}
If $B$ is a compact right-angled polyhedron, then $\Gamma(B)$ is the $1$-skeleton of a $5$-connected triangulation of the $2$-sphere with $m(B)$ vertices.
\end{lemma}

\begin{proof}
By Lemma~\ref{lem:FNS} the triangulation $K(B)$ has no separating triangles. A separating $4$-cycle is either induced, which is excluded by Lemma~\ref{lem:FNS}, or has a chord, and in the latter case one of the two triangles formed with the chord separates. Finally, in a triangulation of the sphere the subgraph induced by a minimal vertex cut contains a cycle separating the two sides of the cut. A cut with at most four vertices would therefore produce a separating cycle of length three or four, which we have just excluded. Since $m(B)\ge12$, the triangulation is $5$-connected.
\end{proof}

By Pogorelov's theorem~\cite{Pog67} the converse holds as well, but we shall not need this: we only use that every facet of $P$ appears in the list of $5$-connected triangulations.

Lemma~\ref{lem:5conn} and Theorem~\ref{thm:andreev} also identify the case of equality in $m(B)\ge12$. If all faces of $B$ are pentagons, then every vertex of the triangulation $K(B)$ has degree five, so $K(B)$ is the icosahedron and $B\cong D$. Thus $D$ is the only compact right-angled polyhedron with twelve faces, and the only one all of whose faces are pentagons.

\subsection{Right-angled polygons}\label{subsec:polygons}
Let $x_1,\dots,x_k$ be the side lengths, in cyclic order, of a compact right-angled polygon. For pentagons and hexagons we have, with indices taken modulo $5$ and $6$ respectively,
\begin{align}
 \cosh x_i&=\sinh x_{i+2}\sinh x_{i+3}=\coth x_{i-1}\coth x_{i+1}&&(k=5),\label{eq:pentagon}\\
 \cosh x_i&=\frac{\cosh x_{i-1}\cosh x_{i+1}+\cosh x_{i+3}}{\sinh x_{i-1}\sinh x_{i+1}}&&(k=6);\label{eq:hexagon}
\end{align}
see, for example,~\cite[Chapter~2]{Bus92}. Let
\[
 \varphi=\frac{1+\sqrt5}{2},\qquad \ell_0=\arcosh\varphi .
\]
The regular right-angled pentagon has side length $\ell_0$, and $\sinh^2\ell_0=\coth^2\ell_0=\varphi$.

\begin{lemma}\label{lem:pentagon}
\begin{enumerate}
\item A right-angled pentagon with two sides of length $\ell_0$ is regular.
\item No compact right-angled polygon with at least six sides has three consecutive sides of length $\ell_0$.
\end{enumerate}
\end{lemma}

\begin{proof}
(i) If the two sides are adjacent, say $x_1=x_2=\ell_0$, the first formula in~\eqref{eq:pentagon} gives $\cosh x_4=\sinh^2\ell_0=\varphi$. If they are not adjacent, say $x_1=x_3=\ell_0$, the second formula gives $\cosh x_2=\coth^2\ell_0=\varphi$. In either case we obtain three sides of length $\ell_0$, two of them adjacent, and repeating the argument shows that all five sides have length $\ell_0$.

(ii) A path of three consecutive sides of length $\ell_0$ meeting at right angles is congruent to three consecutive sides of a regular right-angled pentagon. The lines through the sides immediately before and after this path are then the lines of the remaining two sides of that pentagon, and these lines meet. In a polygon with at least six sides these two sides are not adjacent, so their lines must be ultraparallel by~(F2), a contradiction.
\end{proof}

\begin{lemma}\label{lem:short}
If two consecutive sides of a compact right-angled polygon have lengths $x$ and $y$, then $\sinh x\,\sinh y>1$. The same inequality holds for the lengths of any two edges of $P$ with a common vertex.
\end{lemma}

\begin{proof}
Place the common vertex at $(0,0,1)\in\HH^2$ so that the sides of lengths $x$ and $y$ start in the directions of $-e_2$ and $-e_1$, respectively. The lines containing the next two sides have unit normals $(0,-\cosh x,\sinh x)$ and $(-\cosh y,0,\sinh y)$, whose inner product is $-\sinh x\sinh y$. These lines carry nonadjacent sides of a polygon with at least five sides, so they are ultraparallel and the inner product is less than $-1$. For the second statement, note that two edges of the simple polytope $P$ at a common vertex span a $2$-face, in which they are consecutive sides.
\end{proof}

\subsection{Hyperbolic right-angled L\"obell polyhedra and dodecahedral unions}\label{subsec:unions}
\begin{figure}[ht]
\centering
\begin{tikzpicture}[scale=1.05,line join=round,line cap=round]
\begin{scope}
\def\ra{0.80}\def\rb{1.50}\def\rc{1.95}\def\rd{2.75}
\foreach \i in {0,...,5}{
  \pgfmathsetmacro{\t}{90+60*\i}
  \fill[black!7] (\t:\ra) -- (\t+60:\ra) -- (\t+60:\rb) -- (\t+30:\rc) -- (\t:\rb) -- cycle;
}
\fill[black!22] (90:\ra) \foreach \i in {1,...,5}{ -- (90+60*\i:\ra)} -- cycle;
\foreach \i in {0,...,5}{
  \pgfmathsetmacro{\t}{90+60*\i}
  \draw[line width=1.5pt] (\t:\ra) -- (\t+60:\ra);
  \draw[line width=1.5pt] (\t+30:\rd) -- (\t+90:\rd);
  \draw[line width=0.7pt,densely dashed] (\t:\ra) -- (\t:\rb);
  \draw[line width=0.7pt,densely dashed] (\t+30:\rc) -- (\t+30:\rd);
  \draw[line width=0.5pt] (\t:\rb) -- (\t+30:\rc) -- (\t+60:\rb);
}
\foreach \i in {0,...,5}{
  \pgfmathsetmacro{\t}{90+60*\i}
  \foreach \r in {\ra,\rb}{\fill (\t:\r) circle (0.9pt);}
  \foreach \r in {\rc,\rd}{\fill (\t+30:\r) circle (0.9pt);}
}
\node at (0,-3.25) {(a)};
\end{scope}
\begin{scope}[xshift=5.2cm,yshift=-0.55cm]
\coordinate (p1) at (-0.95,-0.95); \coordinate (p2) at (0.95,-0.95);
\coordinate (p3) at (1.25,0.65);   \coordinate (p4) at (0,1.45);
\coordinate (p5) at (-1.25,0.65);
\fill[black!7] (p1)--(p2)--(p3)--(p4)--(p5)--cycle;
\draw[line width=1.5pt] (p1)--(p2);
\draw[line width=0.7pt,densely dashed] (p2)--(p3) (p5)--(p1);
\draw[line width=0.5pt] (p3)--(p4)--(p5);
\foreach \p in {p1,p2,p3,p4,p5}{\fill (\p) circle (0.9pt);}
\node[below] at ($(p1)!0.5!(p2)$) {$\ell_{n,1}$};
\node[right] at ($(p2)!0.5!(p3)$) {$\ell_{n,2}$};
\node[left]  at ($(p5)!0.5!(p1)$) {$\ell_{n,2}$};
\node[above right] at ($(p3)!0.5!(p4)$) {$\ell_{n,3}$};
\node[above left]  at ($(p4)!0.5!(p5)$) {$\ell_{n,3}$};
\node at (0,-2.7) {(b)};
\end{scope}
\end{tikzpicture}
\caption{(a)~A Schlegel diagram of $L_6$. One base is the outer face, the other is shaded dark, and the two rings of side faces are shaded light and left white. Thick edges are sides of the bases (length $\ell_{6,1}$), dashed edges have exactly one endpoint on a base (length $\ell_{6,2}$), and thin edges are the remaining ones (length $\ell_{6,3}$). (b)~A side face of $L_n$, with the side on a base at the bottom.}\label{fig:lobell}
\end{figure}
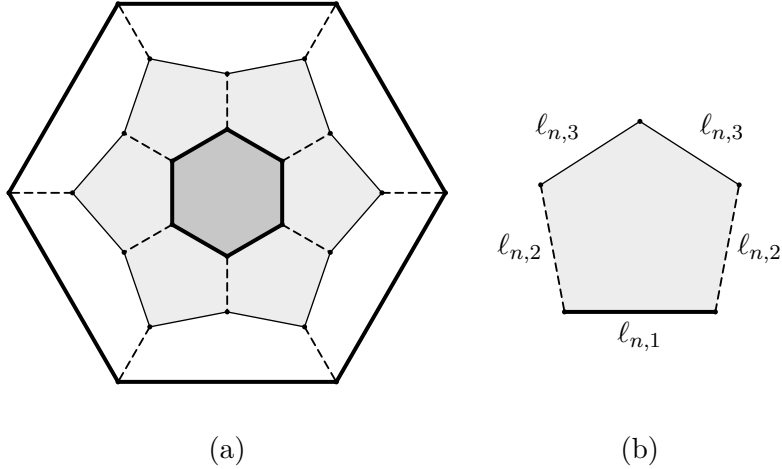

For $n\ge5$, the \emph{L\"obell polyhedron} $L_n$ is the compact right-angled polyhedron with two $n$-gonal faces, its \emph{bases}, and $2n$ pentagonal \emph{side faces}, arranged in two rings of $n$ around the bases~\cite{Loe31,Ves87} (Figure~\ref{fig:lobell}). Thus $m(L_n)=2n+2$, and $L_5=D$ is the right-angled dodecahedron, whose dual graph $\Gamma(D)$ is the graph of the icosahedron. A polyhedron isometric to $D$ is simply called a dodecahedron. The edges of $L_n$ are of three kinds: the sides of the bases, the $2n$ edges with exactly one endpoint on a base, and the remaining $2n$ edges. The symmetries of $L_n$, which are isometries by Theorem~\ref{thm:andreev}, act transitively on the edges of each kind, and we denote the corresponding lengths by $\ell_{n,1}$, $\ell_{n,2}$ and $\ell_{n,3}$. Every side face has one side on a base, and its sides have lengths $(\ell_{n,1},\ell_{n,2},\ell_{n,3},\ell_{n,3},\ell_{n,2})$ in cyclic order, starting from that side.

Two right-angled polytopes can be glued along a common facet. Let $B_1,B_2\subset\HH^n$ be compact right-angled polytopes lying on opposite sides of a hyperplane $\Pi$, with $B_1\cap\Pi=B_2\cap\Pi=F$ a common facet. For each facet $G$ of $F$, the facets of $B_1$ and of $B_2$ through $G$ both lie in the hyperplane through $G$ orthogonal to $\Pi$. Hence $B_1\cup B_2$ is again a compact right-angled polytope. Its facets are those of $B_1$ and $B_2$ other than $F$, the two facets through each facet of $F$ being merged into one.

\begin{figure}[h]
\centering
\begin{tikzpicture}[scale=1.45,line width=0.5pt,line join=round]
\begin{scope}[xshift=0.00cm]
  \fill[black!18] (-0.200,0.275) -- (-0.323,-0.105) -- (-0.000,-0.340) -- (0.323,-0.105) -- (0.200,0.275) -- cycle;
  \draw (0.000,1.000) -- (-0.951,0.309) (0.000,1.000) -- (0.000,0.640) (0.000,0.640) -- (-0.376,0.518) (-0.376,0.518) -- (-0.609,0.198) (-0.376,0.518) -- (-0.200,0.275) (-0.200,0.275) -- (-0.323,-0.105) (-0.951,0.309) -- (-0.588,-0.809) (-0.951,0.309) -- (-0.609,0.198) (-0.609,0.198) -- (-0.609,-0.198) (-0.609,-0.198) -- (-0.376,-0.518) (-0.609,-0.198) -- (-0.323,-0.105) (-0.323,-0.105) -- (-0.000,-0.340) (-0.588,-0.809) -- (0.588,-0.809) (-0.588,-0.809) -- (-0.376,-0.518) (-0.376,-0.518) -- (-0.000,-0.640) (-0.000,-0.640) -- (0.376,-0.518) (-0.000,-0.640) -- (-0.000,-0.340) (-0.000,-0.340) -- (0.323,-0.105) (0.588,-0.809) -- (0.951,0.309) (0.588,-0.809) -- (0.376,-0.518) (0.376,-0.518) -- (0.609,-0.198) (0.609,-0.198) -- (0.609,0.198) (0.609,-0.198) -- (0.323,-0.105) (0.323,-0.105) -- (0.200,0.275) (0.951,0.309) -- (0.000,1.000) (0.951,0.309) -- (0.609,0.198) (0.609,0.198) -- (0.376,0.518) (0.376,0.518) -- (0.000,0.640) (0.376,0.518) -- (0.200,0.275) (0.200,0.275) -- (-0.200,0.275);
\node at (0,-1.2) {$E$};
\end{scope}
\begin{scope}[xshift=2.75cm]
  \fill[black!18] (-0.376,0.518) -- (-0.200,0.275) -- (-0.323,-0.105) -- (-0.609,-0.198) -- (-0.609,0.198) -- cycle;
  \fill[black!18] (-0.588,-0.809) -- (0.588,-0.809) -- (0.376,-0.518) -- (-0.000,-0.640) -- (-0.376,-0.518) -- cycle;
  \draw (0.000,1.000) -- (-0.951,0.309) (0.000,1.000) -- (0.000,0.640) (0.000,0.640) -- (-0.376,0.518) (-0.376,0.518) -- (-0.609,0.198) (-0.376,0.518) -- (-0.200,0.275) (-0.200,0.275) -- (-0.323,-0.105) (-0.951,0.309) -- (-0.588,-0.809) (-0.951,0.309) -- (-0.609,0.198) (-0.609,0.198) -- (-0.609,-0.198) (-0.609,-0.198) -- (-0.376,-0.518) (-0.609,-0.198) -- (-0.323,-0.105) (-0.323,-0.105) -- (-0.000,-0.340) (-0.588,-0.809) -- (0.588,-0.809) (-0.588,-0.809) -- (-0.376,-0.518) (-0.376,-0.518) -- (-0.000,-0.640) (-0.000,-0.640) -- (0.376,-0.518) (-0.000,-0.640) -- (-0.000,-0.340) (-0.000,-0.340) -- (0.323,-0.105) (0.588,-0.809) -- (0.951,0.309) (0.588,-0.809) -- (0.376,-0.518) (0.376,-0.518) -- (0.609,-0.198) (0.609,-0.198) -- (0.609,0.198) (0.609,-0.198) -- (0.323,-0.105) (0.323,-0.105) -- (0.200,0.275) (0.951,0.309) -- (0.000,1.000) (0.951,0.309) -- (0.609,0.198) (0.609,0.198) -- (0.376,0.518) (0.376,0.518) -- (0.000,0.640) (0.376,0.518) -- (0.200,0.275) (0.200,0.275) -- (-0.200,0.275);
\node at (0,-1.2) {$T_2$};
\end{scope}
\begin{scope}[xshift=5.50cm]
  \fill[black!18] (-0.376,0.518) -- (-0.200,0.275) -- (-0.323,-0.105) -- (-0.609,-0.198) -- (-0.609,0.198) -- cycle;
  \fill[black!18] (0.588,-0.809) -- (0.951,0.309) -- (0.609,0.198) -- (0.609,-0.198) -- (0.376,-0.518) -- cycle;
  \draw (0.000,1.000) -- (-0.951,0.309) (0.000,1.000) -- (0.000,0.640) (0.000,0.640) -- (-0.376,0.518) (-0.376,0.518) -- (-0.609,0.198) (-0.376,0.518) -- (-0.200,0.275) (-0.200,0.275) -- (-0.323,-0.105) (-0.951,0.309) -- (-0.588,-0.809) (-0.951,0.309) -- (-0.609,0.198) (-0.609,0.198) -- (-0.609,-0.198) (-0.609,-0.198) -- (-0.376,-0.518) (-0.609,-0.198) -- (-0.323,-0.105) (-0.323,-0.105) -- (-0.000,-0.340) (-0.588,-0.809) -- (0.588,-0.809) (-0.588,-0.809) -- (-0.376,-0.518) (-0.376,-0.518) -- (-0.000,-0.640) (-0.000,-0.640) -- (0.376,-0.518) (-0.000,-0.640) -- (-0.000,-0.340) (-0.000,-0.340) -- (0.323,-0.105) (0.588,-0.809) -- (0.951,0.309) (0.588,-0.809) -- (0.376,-0.518) (0.376,-0.518) -- (0.609,-0.198) (0.609,-0.198) -- (0.609,0.198) (0.609,-0.198) -- (0.323,-0.105) (0.323,-0.105) -- (0.200,0.275) (0.951,0.309) -- (0.000,1.000) (0.951,0.309) -- (0.609,0.198) (0.609,0.198) -- (0.376,0.518) (0.376,0.518) -- (0.000,0.640) (0.376,0.518) -- (0.200,0.275) (0.200,0.275) -- (-0.200,0.275);
\node at (0,-1.2) {$T_3$};
\end{scope}
\end{tikzpicture}
\caption{Schlegel diagrams of $D$ showing the attaching faces. For $E$ one copy of $D$ is shown, with the shaded face glued to a second copy. For $\Tt$ and $\Th$ the central copy is shown, and copies of $D$ are attached along the two shaded faces, which are at distance $2$, respectively $3$, in $\Gamma(D)$.}\label{fig:types}
\end{figure}
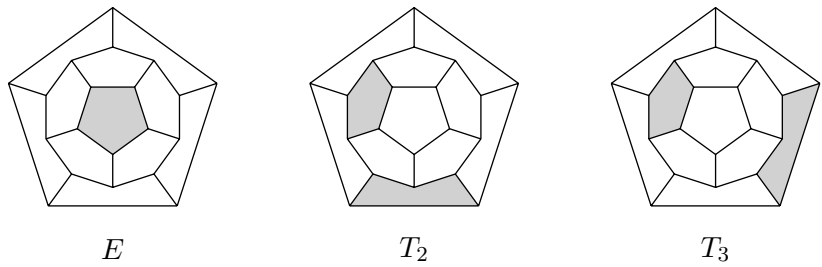

The \emph{double} $E$ of $D$ is the union of $D$ and its reflection in the plane of one of its faces. For the \emph{triple unions}, let $Q_1,Q_2$ be disjoint faces of $D$, with planes $\Pi_1,\Pi_2$ and reflections $\rho_1,\rho_2$ in these planes. The distance between $Q_1$ and $Q_2$ in $\Gamma(D)$ is $2$ or $3$, and accordingly we put
\[
 T_d=\rho_1(D)\cup D\cup\rho_2(D)\qquad(d=2,3).
\]
By~(F2) in Section~\ref{subsec:ra}, the planes $\Pi_1$ and $\Pi_2$ are disjoint, and $\Pi_2$ lies on the same side of $\Pi_1$ as $D$, since it contains $Q_2$. Hence the closed halfspace beyond $\Pi_1$, which contains $\rho_1(D)$ and $Q_1$, does not meet $\Pi_2$ and lies on the same side of $\Pi_2$ as $D$. Thus $D\cup\rho_1(D)$ meets $\Pi_2$ exactly in $Q_2$, we may glue $\rho_2(D)$ to it, and $\Tt$ and $\Th$ are compact right-angled polyhedra. Up to isometry they do not depend on the choice of $Q_1,Q_2$, since the symmetry group of $D$ acts transitively on pairs of faces at a given distance. See Figure~\ref{fig:types}.

In $E$, $\Tt$ and $\Th$, the faces along which the copies of $D$ are glued are called \emph{internal pentagons}. The polyhedron $E$ has one, and each triple union has two, lying in disjoint planes. The plane of an internal pentagon $Q$ meets the polyhedron exactly in $Q$ and is orthogonal to every face that $Q$ meets. Counting faces gives Table~\ref{tab:models}. Every pentagonal face of a member of $\cA=\{D,E,\Tt,\Th\}$ is a regular pentagon, and every hexagonal face is the union of two regular pentagons glued along a side. Every heptagonal face, which occurs only in $\Tt$, is the union of three regular pentagons, the middle one glued to the other two along nonadjacent sides.

\begin{table}[ht]
\centering
\caption{The members of $\cA$.}\label{tab:models}
\begin{tabular}{lcccccc}\toprule
 & copies of $D$ & faces & pentagons & hexagons & heptagons & internal pentagons\\\midrule
$D$  &1&12&12&0&0&0\\
$E$  &2&17&12&5&0&1\\
$\Tt$&3&22&14&6&2&2\\
$\Th$&3&22&12&10&0&2\\\bottomrule
\end{tabular}
\end{table}

\section{A two-layer count in the dual sphere}\label{sec:identity}

\subsection{Defects}
Let $B$ be a compact right-angled polyhedron. Since $B$ is simple, $2f_1(B)=3f_0(B)=\sum_{F}k(F)$, where $F$ ranges over the faces of $B$, and Euler's formula $f_0-f_1+m=2$ yields
\begin{equation}\label{eq:euler}
 m(B)-12=\sum\limits_{F}\bigl(k(F)-5\bigr),
\end{equation}
the sum running over all faces of $B$. All terms are nonnegative, so $m(B)\ge12$, with equality if and only if all faces of $B$ are pentagons, i.e., $B$ is a dodecahedron, by Lemma~\ref{lem:5conn}.

Fix a face $Q$ of $B$ as the \emph{root}, and put $k=k(Q)$. The other faces of $B$ either meet $Q$ along an edge or are disjoint from it. We set
\[
 \dd_1(B,Q)=\sum_{F\sim Q}\bigl(k(F)-5\bigr),\qquad
 \dd_{\ge2}(B,Q)=\sum_{F\cap Q=\varnothing}\bigl(k(F)-5\bigr),
\]
where $F\sim Q$ means that $F$ and $Q$ share an edge. The indices refer to the distance from $Q$ in $\Gamma(B)$. The \emph{defect} of $B$ at $Q$ is defined to be
\begin{equation}\label{eq:defect}
 \dd(B,Q)=m(B)-12+\dd_{\ge2}(B,Q)=(k-5)+\dd_1(B,Q)+2\,\dd_{\ge2}(B,Q),
\end{equation}
the second equality due to~\eqref{eq:euler}. Thus $\dd(B,Q)$ is a nonnegative integer which counts the edge excess $k(F)-5$ of the faces of $B$ over pentagons, with faces far from the root counted twice. By~\eqref{eq:euler} and the remark after Lemma~\ref{lem:5conn}, it vanishes exactly when $B\cong D$. The double weight is not arbitrary but is forced by the count in Proposition~\ref{prop:identity}.

For example, the defect of $E$ is $10$ at the two pentagons disjoint from all hexagons, $8$ at the other ten pentagons, and $7$ at the five hexagons. This is easily read off from Table~\ref{tab:models} and the description of $E$ (see also Appendix~\ref{subsec:small-data}), and will be used repeatedly.

\begin{figure}[h]
\centering
\begin{tikzpicture}[scale=1.25,line width=0.5pt,vtx/.style={circle,fill=black,inner sep=1.1pt}]
  \fill[black!9,even odd rule] (0,0) circle (2.75) (1.146,1.578) -- (0.000,2.120) -- (-1.146,1.578) -- (-1.855,-0.603) -- (-1.585,-1.378) -- (-0.821,-1.933) -- (-0.000,-1.950) -- (1.855,-0.603) -- cycle;
  \fill[black!4] (0.000,1.000) -- (-0.951,0.309) -- (-0.588,-0.809) -- (0.588,-0.809) -- (0.951,0.309) -- cycle;
  \draw (0,0) -- (0.000,1.000);
  \draw[line width=1.1pt] (0.000,1.000) -- (-0.951,0.309);
  \draw (0.000,1.000) -- (-1.146,1.578);
  \draw (0.000,1.000) -- (1.146,1.578);
  \draw (0.000,1.000) -- (0.000,2.120);
  \draw (0,0) -- (-0.951,0.309);
  \draw[line width=1.1pt] (-0.951,0.309) -- (-0.588,-0.809);
  \draw (-0.951,0.309) -- (-1.855,-0.603);
  \draw (-0.951,0.309) -- (-1.146,1.578);
  \draw (0,0) -- (-0.588,-0.809);
  \draw[line width=1.1pt] (-0.588,-0.809) -- (0.588,-0.809);
  \draw (-0.588,-0.809) -- (-0.000,-1.950);
  \draw (-0.588,-0.809) -- (-1.855,-0.603);
  \draw (-0.588,-0.809) -- (-0.821,-1.933);
  \draw (-0.588,-0.809) -- (-1.585,-1.378);
  \draw (0,0) -- (0.588,-0.809);
  \draw[line width=1.1pt] (0.588,-0.809) -- (0.951,0.309);
  \draw (0.588,-0.809) -- (1.855,-0.603);
  \draw (0.588,-0.809) -- (-0.000,-1.950);
  \draw (0,0) -- (0.951,0.309);
  \draw[line width=1.1pt] (0.951,0.309) -- (0.000,1.000);
  \draw (0.951,0.309) -- (1.146,1.578);
  \draw (0.951,0.309) -- (1.855,-0.603);
  \draw[line width=1.1pt,densely dashed] (1.146,1.578) -- (0.000,2.120) -- (-1.146,1.578) -- (-1.855,-0.603) -- (-1.585,-1.378) -- (-0.821,-1.933) -- (-0.000,-1.950) -- (1.855,-0.603) -- cycle;
  \node[vtx,label={[font=\footnotesize]below:$q$}] at (0,0) {};
  \node[vtx] at (0.000,1.000) {};
  \node[vtx] at (-0.951,0.309) {};
  \node[vtx] at (-0.588,-0.809) {};
  \node[vtx] at (0.588,-0.809) {};
  \node[vtx] at (0.951,0.309) {};
  \node[vtx] at (-1.146,1.578) {};
  \node[vtx] at (-1.855,-0.603) {};
  \node[vtx] at (-0.000,-1.950) {};
  \node[vtx] at (1.855,-0.603) {};
  \node[vtx] at (1.146,1.578) {};
  \node[circle,draw,fill=white,inner sep=1.1pt] at (0.000,2.120) {};
  \node[circle,draw,fill=white,inner sep=1.1pt] at (-0.821,-1.933) {};
  \node[circle,draw,fill=white,inner sep=1.1pt] at (-1.585,-1.378) {};
  \node[font=\footnotesize,above right=1pt] at (0.000,1.000) {$u_0$};
  \node[font=\footnotesize,left=2pt] at (-0.951,0.309) {$u_1$};
  \node[font=\footnotesize] at (-0.330,-0.990) {$u_2$};
  \node[font=\footnotesize] at (0.330,-0.990) {$u_3$};
  \node[font=\footnotesize,right=2pt] at (0.951,0.309) {$u_{4}$};
  \node[font=\footnotesize,above left=1pt] at (-1.146,1.578) {$x_0$};
  \node[font=\footnotesize,left=2pt] at (-1.855,-0.603) {$x_1$};
  \node[font=\footnotesize,below=2pt] at (-0.000,-1.950) {$x_2$};
  \node[font=\footnotesize,right=2pt] at (1.855,-0.603) {$x_3$};
  \node[font=\footnotesize,above right=1pt] at (1.146,1.578) {$x_4$};
  \node[font=\footnotesize] at (1.150,-2.300) {$\mathcal{D}$};
\end{tikzpicture}
\caption{The collar of a root $q$ in $K(B)$, here with $k=5$. The bold cycle is $C=\lk(q)$, and the dashed cycle is the outer path $\gamma$. The vertex $x_i$ is adjacent to $u_i$ and $u_{i+1}$. The open vertices on $\gamma$ between $x_{i-1}$ and $x_i$ are adjacent to $u_i$ only, and there are $k(u_i)-5$ of them. The shaded region $\mathcal D$ is the complementary disk, whose interior vertices are the $n_{\ge3}$ vertices at distance at least three from $q$.}\label{fig:collar}
\end{figure}
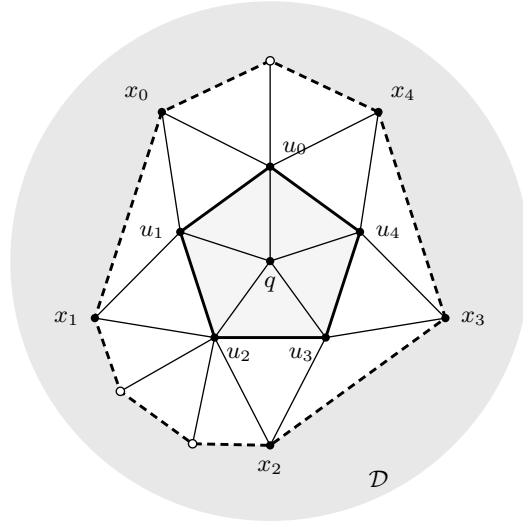

\subsection{The collar of a face}
The following lemma describes the faces of $B$ near a root. In Section~\ref{subsec:identity} it will be applied to every neighbour of a facet of $P$.

\begin{lemma}\label{lem:collar}
Let $Q$ be a $k$-gonal face of a compact right-angled hyperbolic polyhedron $B$, and let $n_{\ge3}(B,Q)$ be the number of faces at distance at least three from $Q$ in $\Gamma(B)$. Then
\begin{equation}\label{eq:collar}
 m(B)=1+2k+\dd_1(B,Q)+n_{\ge3}(B,Q),\qquad \dd_{\ge2}(B,Q)=n_{\ge3}(B,Q)+k-6.
\end{equation}
Moreover $n_{\ge3}(B,Q)\ge1$, and equality holds only when $\dd_1(B,Q)=0$ and $B$ is isometric to the L\"obell polyhedron $L_k$. Consequently, if $B\not\cong L_k$, then
\begin{equation}\label{eq:collar-bound}
 \dd_{\ge2}(B,Q)\ge k-4\qquad\text{and}\qquad m(B)-12\ge 2k-9 .
\end{equation}
In particular $\dd_{\ge2}(B,Q)\ge1$ whenever $B\not\cong D$.
\end{lemma}

\begin{proof}
Let $q$ be the vertex dual to $Q$. Then $\lk(q)=(u_0,\dots,u_{k-1})$ is an induced cycle, with indices modulo $k$, as shown in Figure~\ref{fig:collar}.
Let $x$ be a vertex not in $\{q\}\cup \lk(q)$. If $x$ had two nonadjacent neighbours $u_i,u_j$ on $\lk(q)$, then $x,u_i,q,u_j$ would be an induced square. Hence the neighbours of $x$ on $\lk(q)$ form a clique of the induced cycle $\lk(q)$, so $x$ has at most two of them, and if two, they are consecutive. 

For each $i$ the edge $u_iu_{i+1}$ lies in exactly two triangles, say $qu_iu_{i+1}$ and $x_iu_iu_{i+1}$. The vertices $x_i$ are not on $\lk(q)$, and by flagness they are exactly the vertices with two neighbours on $\lk(q)$.

The link of $u_i$ is a cycle of length $k(u_i)$ containing $q,u_{i-1},u_{i+1},x_{i-1},x_i$. Its remaining $k(u_i)-5$ vertices are adjacent to $u_i$ and, by the above, to no other vertex of $\lk(q)$. Conversely, every vertex with $u_i$ as its only neighbour on $\lk(q)$ lies in $\lk(u_i)$. So exactly $\dd_1(B,Q)=\sum_i(k(u_i)-5)$ vertices have one neighbour on $\lk(q)$. Counting $q$, the $k$ vertices of $\lk(q)$, the $k$ vertices $x_i$, these $\dd_1(B,Q)$ vertices and the $n_{\ge3}(B,Q)$ vertices with no neighbour on $\lk(q)$ gives the first identity in~\eqref{eq:collar}. The second follows from the first and~\eqref{eq:defect}.

In $\lk(u_i)$, the arc from $x_{i-1}$ to $x_i$ that avoids $q$ has as interior vertices precisely the vertices whose only neighbour on $\lk(q)$ is $u_i$, as shown in Figure~\ref{fig:collar}. The concatenation of these arcs is a simple closed curve $\gamma$ in $K(B)$ of length $k+\dd_1(B,Q)$. The union of the closed stars of $q,u_0,\dots,u_{k-1}$ is a disk $\mathcal{S}$ with boundary $\gamma$. Let $\mathcal D$ be the complementary closed disk, which is also bounded by $\gamma$. Interior vertices of $\mathcal{D}$ are exactly the $n_{\ge3}(B,Q)$ vertices at distance at least three from $q$.

Since every vertex of $K(B)$ has degree at least five, the vertices $x_i$ and the other vertices of $\gamma$ have degree at least one and two in the interior of $\mathcal D$, respectively.
This is impossible when $n_{\ge3}(B,Q)=0$, i.e., when $\mathcal D$ is a triangulated polygon, and no interior vertex can be found. If $n_{\ge3}(B,Q)=1$, we only have one vertex $z$ in the interior of $\mathcal{D}$, and only the vertices $x_i$ survive, all with degree five. The complex $K(B)$ now consists of the two cones over the $k$-cycles $\lk(q)$ and $\gamma$ with apices $q$ and $z$, joined by an antiprism. This is $K(L_k)$, and $B\cong L_k$ by Theorem~\ref{thm:andreev}.

If $B\not\cong L_k$ we therefore have $n_{\ge3}(B,Q)\ge2$, and~\eqref{eq:collar-bound} follows from~\eqref{eq:collar} and~\eqref{eq:defect}.
\end{proof}

\subsection{The defect identity}\label{subsec:identity}
We now return to the $4$-polytope $P$ and write $K=K(P)$. For a facet $A$ of $P$ we define the \emph{defect of $A$ in $P$} by
\[
 \dd_P(A)=\sum_{B\sim A}\dd(B,A\cap B),
\]
the sum running over the neighbours of $A$.

\begin{proposition}\label{prop:identity}
Let $A$ be a facet of $P$ and let $n_{\ge3}(A,P)$ be the number of facets of $P$ at distance at least three from $A$ in $\Gamma(P)$. Then
\begin{equation}\label{eq:identity}
 2N=30+5\,m(A)+\dd_P(A)+2\,n_{\ge3}(A).
\end{equation}
\end{proposition}

\begin{proof}
Let $v$ be the vertex of $K$ dual to $A$, so that $\lk_K(v)\cong K(A)$. Its $m=m(A)$ vertices are dual to the neighbours of $A$. Call the remaining vertices of $K$, other than $v$ and those of $\lk_K(v)$, \emph{outer vertices}. The neighbours in $\lk_K(v)$ of an outer vertex $x$ are pairwise adjacent, since two nonadjacent ones would form an induced square with $x$ and $v$. As $K$ is flag and $3$-dimensional, $x$ has at most three neighbours in $\lk_K(v)$, and they span a simplex. For $0\le j\le3$ let $V_j$ be the set of outer vertices with exactly $j$ neighbours in $\lk_K(v)$. The vertices in $V_0$ are those at distance at least three from $v$, so $|V_0|=n_{\ge3}(A,P)$, and
\begin{equation}\label{eq:N-count}
 N=1+m+|V_1|+|V_2|+|V_3|+n_{\ge3}(A,P) .
\end{equation}

Let $u\in \lk_K(v)$ be dual to the neighbour $B$ of $A$, and let $Q=A\cap B$, so that $\lk_K(u)\cong K(B)$ and $v$ is the vertex of $\lk_K(u)$ dual to the face $Q$ of $B$. The vertices of $\lk_K(u)$ adjacent to $v$ form the cycle $C=\lk_K(u)\cap \lk_K(v)$, and all other vertices of $\lk_K(u)$ except $v$ are outer. By the proof of Lemma~\ref{lem:collar}, applied to $B$ with root $Q$, these outer vertices are of three kinds: $k(Q)$ vertices adjacent to two vertices of $C$, $\dd_1(B,Q)$ vertices adjacent to one vertex of $C$, and $n_{\ge3}(B,Q)$ vertices adjacent to no vertex of $C$. If $x$ is an outer vertex adjacent to $u$, then its neighbours in $\lk_K(v)$ are pairwise adjacent, so those other than $u$ lie in $\lk_K(u)\cap \lk_K(v)=C$. Hence the neighbours of $x$ in $\lk_K(v)$ are $u$ and the neighbours of $x$ on $C$. Thus the three kinds belong to $V_3$, $V_2$ and $V_1$, respectively. Conversely, a vertex in $V_j$ with $j\ge1$ lies in the links of exactly $j$ vertices of $\lk_K(v)$. Summing over the $m$ neighbours $B$ of $A$, we obtain
\[
 3|V_3|=\sum_{B\sim A}k(A\cap B),\qquad 2|V_2|=\sum_{B\sim A}\dd_1(B,A\cap B),\qquad |V_1|=\sum_{B\sim A}n_{\ge3}(B,A\cap B).
\]
The first sum is $\sum_Fk(F)=6m-12$, the sum over the faces of $A$, so $|V_3|=2m-4$. By Lemma~\ref{lem:collar}, $n_{\ge3}(B,Q)=\dd_{\ge2}(B,Q)-k(Q)+6$, so $|V_1|=\sum_{B\sim A}\dd_{\ge2}(B,A\cap B)+12$. Substituting into~\eqref{eq:N-count},
\[
 2N=18+6m+2\,n_{\ge3}(A)+\sum_{B\sim A}\bigl(\dd_1(B,A\cap B)+2\,\dd_{\ge2}(B,A\cap B)\bigr).
\]
By~\eqref{eq:defect} the last sum is $\dd_P(A)-\sum_F(k(F)-5)$, and $\sum_F(k(F)-5)=m-12$ by~\eqref{eq:euler}. This gives~\eqref{eq:identity}.
\end{proof}

For the right-angled $120$-cell, $m(A)=12$ and $\dd_Z(A)=0$ for every facet, and~\eqref{eq:identity} says that $75$ of the $120$ facets lie at distance at least three from any given one. In general, since $n_{\ge3}(A)\ge0$, Proposition~\ref{prop:identity} gives:

\begin{corollary}\label{cor:budget}
If $N\le120$, then every facet $A$ of a compact right-angled hyperbolic 4-polytope $P$ satisfies
$
 \dd_P(A)\le210-5\,m(A).
$
\end{corollary}

\section{Facets with fewer than twenty faces}\label{sec:small}

The aim of this section is to prove the following proposition.

\begin{proposition}\label{prop:small}
Suppose $P$ is a compact right-angled hyperbolic $4$-polytope with $N\le120$ facets. Then every facet of $P$ is isometric to the right-angled dodecahedron $D$, to its double $E$, or to a polyhedron with at least twenty faces. Moreover, no facet of $P$ is isometric to a L\"obell polyhedron $L_n$ with $n\ge6$.
\end{proposition}

Although a single face of a small polyhedron $A$ may well admit a neighbour of small defect, several such neighbours are mutually exclusive. We need to control both individual neighbours and their compatibility.

\subsection{The small polyhedra}\label{subsec:small-list}
By Lemma~\ref{lem:5conn} the combinatorial types of compact right-angled polyhedra with $m$ faces are among the $5$-connected triangulations of the sphere with $m$ vertices. For $12\le m\le19$ there are
\[
 1,\ 0,\ 1,\ 1,\ 3,\ 4,\ 12,\ 23
\]
of them~\cite{BM07} (see Appendix~\ref{subsec:catalogues}). Those with at most $17$ faces are listed in Table~\ref{tab:small}, and their rotation systems are reproduced in Appendix~\ref{subsec:small-data}. Besides the dodecahedron $D$, its double $E$ and the L\"obell polyhedra $L_6$ and $L_7$, there is one type $F_{15}$ with $15$ faces, two types $F_{16a}$, $F_{16b}$ with $16$ faces, and three types $F_{17a},F_{17b},F_{17c}$ with $17$ faces. Note that $F_{16b}$ has twelve pentagons, four pairwise disjoint hexagons, and every pentagon of $F_{16b}$ meets exactly two hexagons.
The label $(m,i)$ in the table is the position of the type in the catalogue described in Appendix~\ref{subsec:catalogues}, and we use these labels when referring to the data.

\begin{table}[ht]
\centering
\caption{The combinatorial types with at most $17$ faces, all of which admit a compact right-angled hyperbolic structure. The column $(p_5,p_6,p_7)$ gives the numbers of pentagons, hexagons and heptagons, and $|\Aut|$ is the order of the automorphism group of $\Gamma$, which by Theorem~\ref{thm:andreev} is the isometry group. The names in the first column are used throughout the paper.}\label{tab:small}
\begin{tabular}{lcccc}\toprule
 & label $(m,i)$ & faces & $(p_5,p_6,p_7)$ & $|\Aut|$\\\midrule
$D$       &$(12,0)$&12&$(12,0,0)$&120\\
$L_6$     &$(14,0)$&14&$(12,2,0)$&24\\
$F_{15}$  &$(15,0)$&15&$(12,3,0)$&12\\
$L_7$     &$(16,0)$&16&$(14,0,2)$&28\\
$F_{16a}$ &$(16,1)$&16&$(12,4,0)$&4\\
$F_{16b}$ &$(16,2)$&16&$(12,4,0)$&24\\
$F_{17a}$ &$(17,0)$&17&$(13,3,1)$&2\\
$F_{17b}$ &$(17,1)$&17&$(12,5,0)$&4\\
$F_{17c}$ &$(17,2)$&17&$(12,5,0)$&4\\
$E$       &$(17,3)$&17&$(12,5,0)$&20\\\bottomrule
\end{tabular}
\end{table}

\begin{lemma}\label{lem:small-facts}
\begin{enumerate}
\item Among compact right-angled hyperbolic polyhedra $B$ with $14\le m\le19$ faces, the double $E$ of the right-angled dodecahedron is the only one that has a regular pentagonal face.
\item If $m(B)=17$, then $\dd(B,Q)\ge7$ for every face $Q$ of the polyhedron $B$.
\end{enumerate}
\end{lemma}

\begin{proof}
 Part~(ii) is a direct check on the four combinatorial types with $17$ faces (Table~\ref{tab:small}): by~\eqref{eq:defect}, $\dd(B,Q)=5+\dd_{\ge2}(B,Q)$, so one needs $\dd_{\ge2}(B,Q)\ge2$, and this holds at every face of each of them (the rotation systems and the defects at all faces are given in Appendix~\ref{subsec:small-data}). Part~(i) is obtained by a computation in ball arithmetic, as described in Appendix~\ref{app:computation}. For each of the $43$ types with $14\le m\le19$ faces other than $E$, every pentagonal face has a side of length not equal to $\ell_0$. 
\end{proof}

\subsection{Attaching a polyhedron along a face}\label{subsec:placements}

\begin{definition}\label{def:placement}
Let $Q$ be a face of a polyhedron $A$ and $R$ a face of a polyhedron $B$ with $k(R)=k(Q)$. A \emph{placement} of $B$ on $Q$ along $R$ is a combinatorial isomorphism $\iota\colon Q\to R$. There are $2k(Q)$ placements for each pair $(Q,R)$. We call $A$ the \emph{source} and $B$ the \emph{target} of the placement.
\end{definition}

A placement is a purely combinatorial object: it only records which face of $B$ is glued to $Q$ and how the two boundaries are matched, and for given $A$, $Q$ and $B$ there are finitely many placements. A pair of neighbouring facets $A,B$ of a compact right-angled hyperbolic $4$-polytope $P$ gives a placement of $B$ on $A\cap B$, which we call \emph{realized}. Whether a placement can be realized is a matter of metrics.

\begin{figure}[ht]
\centering
\begin{tikzpicture}[scale=1.05,line width=0.5pt,line join=round,line cap=round,
  vx/.style={circle,fill=black,inner sep=1.4pt},
  lab/.style={font=\footnotesize}]
\begin{scope}
  \coordinate (z) at (0,0);
  \coordinate (ee) at (-2.3,0);      
  \coordinate (ff) at (1.05,0.75);   
  \coordinate (dd) at (0,-1.75);     
  \coordinate (tt) at (0,1.75);      
  \fill[black!20] (z) -- (ee) -- ($(ee)+(ff)$) -- (ff) -- cycle;
  \fill[black!7]  (z) -- (ee) -- ($(ee)+(dd)$) -- (dd) -- cycle;
  \fill[black!12] (z) -- (ff) -- ($(ff)+(dd)$) -- (dd) -- cycle;
  \draw[black!45,densely dotted] (ee) -- ($(ee)+(ff)$) -- (ff) -- ($(ff)+(dd)$) -- (dd) -- ($(ee)+(dd)$) -- (ee);
  \draw[line width=1.1pt] (z) -- node[lab,above,pos=0.55] {$e$} (ee);
  \draw[line width=1.1pt] (z) -- (ff);
  \draw[line width=1.1pt] (z) -- node[lab,left,pos=0.6] {$t_A(z)$} (dd);
  \draw[line width=1.6pt,densely dashed] (z) -- node[lab,right,pos=0.6] {$t$} (tt);
  \node[vx] at (tt) {};
  \node[vx,label={[lab]right:$z$}] at (z) {};
  \node[lab] at (-0.75,0.38) {$Q$};
  \node[lab] at (0.55,-0.55) {$Q'$};
  \node[lab] at (-1.55,-0.65) {$A$};
  \node[lab] at (-1.15,1.15) {$B$};
  \node[lab] at (1.05,1.45) {$C$};
  \node[font=\small] at (-0.6,-2.45) {(a)};
\end{scope}
\begin{scope}[xshift=5.6cm]
  \coordinate (z) at (-1.15,0);
  \coordinate (zp) at (1.15,0);
  \coordinate (ff) at (1.3,0.7);
  \coordinate (dd) at (0,-1.45);
  \coordinate (p5) at (-1.15,1.25);
  \coordinate (p3) at (1.15,1.25);
  \coordinate (p4) at (0,2.0);
  \fill[black!20] (z) -- (zp) -- ($(zp)+(ff)$) -- ($(z)+(ff)$) -- cycle;
  \fill[black!7]  (z) -- (zp) -- ($(zp)+(dd)$) -- ($(z)+(dd)$) -- cycle;
  \draw[black!45,densely dotted] (zp) -- ($(zp)+(ff)$) -- ($(z)+(ff)$) -- (z);
  \draw[black!45,densely dotted] ($(z)+(dd)$) -- ($(zp)+(dd)$);
  \fill[black!12,opacity=0.85] (z) -- (zp) -- (p3) -- (p4) -- (p5) -- cycle;
  \draw (p3) -- (p4) -- (p5);
  \draw[line width=1.1pt] (z) -- node[lab,below,pos=0.5] {$e$} (zp);
  \draw[line width=1.6pt,densely dashed] (z) -- node[lab,left,pos=0.75] {$t_B(z)$} (p5);
  \draw[line width=1.6pt,densely dashed] (zp) -- node[lab,right,pos=0.85] {$t_B(z')$} (p3);
  \draw[line width=1.1pt] (z) -- node[lab,left] {$t_A(z)$} ($(z)+(dd)$);
  \draw[line width=1.1pt] (zp) -- node[lab,right] {$t_A(z')$} ($(zp)+(dd)$);
  \node[vx,label={[lab]left:$z$}] at (z) {};
  \node[vx,label={[lab]right:$z'$}] at (zp) {};
  \node[lab] at (0,0.95) {$F'$};
  \node[lab] at (1.6,0.36) {$Q$};
  \node[lab] at (0,-0.85) {$A$};
  \node[font=\small] at (0,-2.45) {(b)};
\end{scope}
\end{tikzpicture}
\caption{The local structure at a vertex $z$ of $P$ lying on $Q=A\cap B$ (Lemma~\ref{lem:placement}), drawn schematically with $A$ below $Q$ and its neighbours above. (a)~Near $z$ the facet $A$ looks like the corner of a cube: it has three faces at $z$, among them $Q$ and $Q'=A\cap C$, and three edges at $z$, namely the two edges of $Q$ and $t_A(z)$. The fourth edge $t$ of $P$ at $z$ (dashed) does not lie in $A$, so it lies in every neighbour of $A$ through $z$, and $t=t_B(z)=t_C(z)$. This is~(iv). Any two of the four edges at $z$ span a $2$-face, in which they are consecutive sides. This gives~(ii). (b)~The face $F'$ of $B$ through an edge $e=zz'$ of $Q$. Its sides at $z$ and $z'$ are $t_B(z)$ and $t_B(z')$, which gives~(iii) when $F'$ is a pentagon.}\label{fig:vertex}
\end{figure}
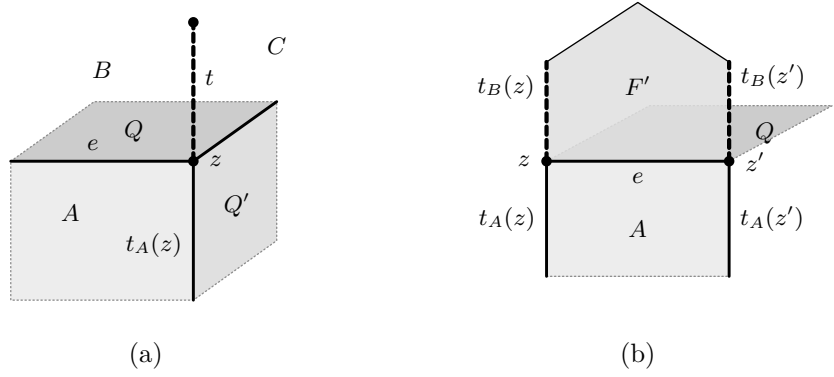

Let \(z\) be a vertex of the polytope \(P\) that lies on \(Q=A\cap B\). Four edges of \(P\) meet at \(z\): two are the edges of the polygon \(Q\) incident to \(z\), and the remaining two are \(t_A(z)\) and \(t_B(z)\), which lie in \(A\) and \(B\), respectively, as shown in Figure~\ref{fig:vertex}. Lemma~\ref{lem:placement} lists the necessary, but not sufficient, conditions that we use to exclude unrealized placements.

\begin{lemma}\label{lem:placement}
Let $A,B$ be neighbouring facets of a compact right-angled hyperbolic $4$-polytope $P$ with $Q=A\cap B$.
\begin{enumerate}
\item Each edge of $Q$ has the same length in $A$ and in $B$.
\item If $z$ is a vertex of $Q$ and $a$ is an edge of $A$ at $z$, then
\[
 \sinh\ell(a)\,\sinh\ell(t_B(z))>1 .
\]
\item Let $e=zz'$ be an edge of $Q$ and $F'$ the face of $B$ through $e$ other than $Q$. If $F'$ is a pentagon, then
\[
 \cosh\ell(e)<\cosh\ell(a)\,\cosh\ell(a')
\]
for every edge $a$ of $A$ at $z$ and every edge $a'$ of $A$ at $z'$.
\item If $C$ is the neighbour of $A$ across another face $Q'$ of $A$ through $z$, then $t_C(z)=t_B(z)$. In particular, the edges $t_B(z)$ and $t_C(z)$ have the same length.
\end{enumerate}
\end{lemma}

\begin{proof}
Part~(i) is clear, and~(ii) is Lemma~\ref{lem:short} applied to $a$ and $t_B(z)$. For~(iii), the sides of the pentagon $F'$ adjacent to $e$ are $t_B(z)$ and $t_B(z')$, so~\eqref{eq:pentagon} gives $\cosh\ell(e)=\coth\ell(t_B(z))\coth\ell(t_B(z'))$. By~(ii),
\[
 \coth^2\ell(t_B(z))=1+\frac{1}{\sinh^2\ell(t_B(z))}<1+\sinh^2\ell(a)=\cosh^2\ell(a),
\]
and similarly at $z'$. Part~(iv) holds because $P$ has only one edge at $z$ outside $A$.
\end{proof}

Conditions~(i)--(iii) concern a single neighbour, and condition~(iv) relates neighbours across adjacent faces of $A$ (Figure~\ref{fig:vertex}). Since $A$ and $B$ are known explicitly up to isometry, each condition can be checked for a given placement. Appendix~\ref{subsec:matching} explains how this is done rigorously.

\subsection{Bounds for neighbour defects}\label{subsec:single}
The first lemma bounds the defect of a single neighbour.

\begin{lemma}\label{lem:single}
Let $A,B$ be neighbouring facets of a compact right-angled hyperbolic $4$-polytope $P$ with $Q=A\cap B$.
\begin{enumerate}
\item If $A\cong L_6$ and $Q$ is a pentagon, then $\dd(B,Q)\ge10$. If $A\cong L_7$ and $Q$ is a pentagon, then $\dd(B,Q)\ge8$.
\item Let $A\cong F_{15}$, $F_{16a}$ or $F_{16b}$. If $B\not\cong A$, then $\dd(B,Q)\ge13$, $10$ or $10$, respectively. If $B\cong A$, then $\dd(B,Q)\ge4$, $6$ or $6$, respectively. Moreover, if $A\cong B\cong F_{16b}$, then $\dd(B,Q)=6$ when $Q$ is a pentagon and $\dd(B,Q)=7$ when $Q$ is a hexagon.
\item If $A\cong F_{17a}$, $F_{17b}$ or $F_{17c}$, then $\dd(B,Q)\ge7$, and equality implies $B\cong A$. A neighbour of defect $7$ can occur across at most six faces of $F_{17a}$, three faces of $F_{17b}$ and two faces of $F_{17c}$, and across all other faces $\dd(B,Q)\ge8$.
\end{enumerate}
\end{lemma}

\begin{proof}
We first reduce each statement to a finite check. If $B\not\cong D$, then $\dd_{\ge2}(B,Q)\ge1$ by Lemma~\ref{lem:collar}, so $\dd(B,Q)\ge m(B)-11$ by~\eqref{eq:defect}. Hence, for any integer $\dd_*$,
\begin{equation}\label{eq:finite}
 \dd(B,Q)<\dd_*\quad\Longrightarrow\quad B\cong D\quad\text{or}\quad m(B)\le\dd_*+10 .
\end{equation}
To prove a bound $\dd(B,Q)\ge\dd_*$ it therefore suffices to consider the finitely many types $B$ with at most $\dd_*+10$ faces, and for each of them the faces $R$ with $k(R)=k(Q)$ and $\dd(B,R)<\dd_*$, together with all $2k(Q)$ placements along $R$. A bound is established by showing that each of these placements violates one of the conditions~(i)--(iii) of Lemma~\ref{lem:placement}. A placement that violates none of them is said to \emph{survive}.

For~(i) the cutoffs $\dd_*=10$ and $8$ lead to targets with at most $20$ and $18$ faces, respectively, and every placement is excluded. For the first assertion of~(ii), the cutoffs $13,10,10$ lead to targets with at most $23,20,20$ faces, and again every placement on a target not isometric to $A$ is excluded. The placements of $A$ on itself that survive all three conditions are listed in Table~\ref{tab:same-type} of Appendix~\ref{app:data}, together with their defects, and this gives the remaining assertions of~(ii). For~(iii) it suffices to examine targets with at most $18$ faces and the cutoff $8$. No placement of defect less than $7$ survives. Every surviving placement of defect $7$ is a placement of $A$ on itself, and the faces of $A$ that carry one are six, three and two in number for $F_{17a},F_{17b},F_{17c}$. The numbers of placements involved in~(i) and~(ii) are given in Table~\ref{tab:coverage}.
\end{proof}

Condition~(iv) of Lemma~\ref{lem:placement} limits the low-defect neighbours that can occur together.

\begin{lemma}\label{lem:simultaneous}
Let $A$ be a facet of $P$.
\begin{enumerate}
\item If $A\cong F_{15}$, then at most six neighbours of $A$ are isometric to $A$.
\item If $A\cong F_{16a}$, then at most four neighbours of $A$ are isometric to $A$.
\item Suppose $A\cong F_{16b}$. Let $Q$ be a pentagon and $R$ a hexagon of $A$ with a common edge, and let $B,C$ be the neighbours of $A$ across $Q$ and $R$. If $B\cong A$, then $C\not\cong A$, and $B\cap C$ is a hexagon.
\item Suppose $A\cong F_{16b}$, and let $C$ be the neighbour across a hexagon $R$ of $A$. If $h$ of the six pentagons of $A$ adjacent to $R$ have neighbours isometric to $A$, then $\dd(C,R)\ge2h+3$.
\end{enumerate}
\end{lemma}

\begin{figure}[ht]
\centering
\begin{tikzpicture}[scale=0.92,line width=0.5pt,line join=round,
  lab/.style={font=\footnotesize}]
\begin{scope}
  \fill[black!5] (0.823,0.475) -- (0.000,0.950) -- (-0.823,0.475) -- (-0.823,-0.475) -- (-0.000,-0.950) -- (0.823,-0.475) -- cycle;
  \fill[black!22] (0.823,0.475) -- (0.000,0.950) -- (0.000,2.050) -- (1.200,2.078) -- (1.775,1.025) -- cycle;
  \fill[black!22] (0.000,0.950) -- (-0.823,0.475) -- (-1.775,1.025) -- (-1.200,2.078) -- (0.000,2.050) -- cycle;
  \fill[black!0] (-0.823,0.475) -- (-0.823,-0.475) -- (-1.775,-1.025) -- (-2.400,0.000) -- (-1.775,1.025) -- cycle;
  \fill[black!22] (-0.823,-0.475) -- (-0.000,-0.950) -- (-0.000,-2.050) -- (-1.200,-2.078) -- (-1.775,-1.025) -- cycle;
  \fill[black!0] (-0.000,-0.950) -- (0.823,-0.475) -- (1.775,-1.025) -- (1.200,-2.078) -- (-0.000,-2.050) -- cycle;
  \fill[black!0] (0.823,-0.475) -- (0.823,0.475) -- (1.775,1.025) -- (2.400,-0.000) -- (1.775,-1.025) -- cycle;
  \draw (0.000,0.950) -- (0.000,2.050) -- (1.200,2.078) -- (1.775,1.025);
  \draw (-0.823,0.475) -- (-1.775,1.025) -- (-1.200,2.078) -- (0.000,2.050);
  \draw (-0.823,-0.475) -- (-1.775,-1.025) -- (-2.400,0.000) -- (-1.775,1.025);
  \draw (-0.000,-0.950) -- (-0.000,-2.050) -- (-1.200,-2.078) -- (-1.775,-1.025);
  \draw (0.823,-0.475) -- (1.775,-1.025) -- (1.200,-2.078) -- (-0.000,-2.050);
  \draw (0.823,0.475) -- (1.775,1.025) -- (2.400,-0.000) -- (1.775,-1.025);
  \draw[line width=1.0pt] (0.823,0.475) -- (0.000,0.950) -- (-0.823,0.475) -- (-0.823,-0.475) -- (-0.000,-0.950) -- (0.823,-0.475) -- cycle;
  \draw (0.823,0.475) -- (1.775,1.025);
  \draw (0.000,0.950) -- (0.000,2.050);
  \draw (-0.823,0.475) -- (-1.775,1.025);
  \draw (-0.823,-0.475) -- (-1.775,-1.025);
  \draw (-0.000,-0.950) -- (-0.000,-2.050);
  \draw (0.823,-0.475) -- (1.775,-1.025);
  \node[lab] at (0,0) {$R$};
  \node[lab] at (0.810,1.403) {$Q_1$};
  \node[lab] at (-0.810,1.403) {$Q_2$};
  \node[lab] at (-1.620,0.000) {$Q_3$};
  \node[lab] at (-0.810,-1.403) {$Q_4$};
  \node[lab] at (0.810,-1.403) {$Q_5$};
  \node[lab] at (1.620,-0.000) {$Q_6$};
  \node[font=\small] at (0,-2.95) {(a)~~in $A$};
\end{scope}
\begin{scope}[xshift=6.3cm]
  \fill[black!5] (0.823,0.475) -- (0.000,0.950) -- (-0.823,0.475) -- (-0.823,-0.475) -- (-0.000,-0.950) -- (0.823,-0.475) -- cycle;
  \fill[black!22] (0.823,0.475) -- (0.000,0.950) -- (0.000,2.050) -- (0.702,2.295) -- (1.637,1.755) -- (1.775,1.025) -- cycle;
  \fill[black!22] (0.000,0.950) -- (-0.823,0.475) -- (-1.775,1.025) -- (-1.637,1.755) -- (-0.702,2.295) -- (0.000,2.050) -- cycle;
  \fill[black!22] (-0.823,-0.475) -- (-0.000,-0.950) -- (-0.000,-2.050) -- (-0.702,-2.295) -- (-1.637,-1.755) -- (-1.775,-1.025) -- cycle;
  \draw (0.000,0.950) -- (0.000,2.050) -- (0.702,2.295) -- (1.637,1.755) -- (1.775,1.025);
  \draw (-0.823,0.475) -- (-1.775,1.025) -- (-1.637,1.755) -- (-0.702,2.295) -- (0.000,2.050);
  \draw[densely dotted] (-1.775,-1.025) -- (-2.400,0.000) -- (-1.775,1.025);
  \draw (-0.000,-0.950) -- (-0.000,-2.050) -- (-0.702,-2.295) -- (-1.637,-1.755) -- (-1.775,-1.025);
  \draw[densely dotted] (1.775,-1.025) -- (1.200,-2.078) -- (-0.000,-2.050);
  \draw[densely dotted] (1.775,1.025) -- (2.400,-0.000) -- (1.775,-1.025);
  \draw[line width=1.0pt] (0.823,0.475) -- (0.000,0.950) -- (-0.823,0.475) -- (-0.823,-0.475) -- (-0.000,-0.950) -- (0.823,-0.475) -- cycle;
  \draw (0.823,0.475) -- (1.775,1.025);
  \draw (0.000,0.950) -- (0.000,2.050);
  \draw (-0.823,0.475) -- (-1.775,1.025);
  \draw (-0.823,-0.475) -- (-1.775,-1.025);
  \draw (-0.000,-0.950) -- (-0.000,-2.050);
  \draw (0.823,-0.475) -- (1.775,-1.025);
  \node[lab] at (0,0) {$R$};
  \node[lab] at (0.830,1.438) {$B_1\cap C$};
  \node[lab] at (-0.830,1.438) {$B_2\cap C$};
  \node[lab] at (-1.620,0.000) {?};
  \node[lab] at (-0.830,-1.438) {$B_4\cap C$};
  \node[lab] at (0.810,-1.403) {?};
  \node[lab] at (1.620,-0.000) {?};
  \node[font=\small] at (0,-2.95) {(b)~~in $C$};
\end{scope}
\end{tikzpicture}
\caption{Lemma~\ref{lem:simultaneous}(iv), drawn schematically for $h=3$. (a)~A hexagon $R$ of the facet $A\cong F_{16b}$ with the six pentagons $Q_1,\dots,Q_6$ of $A$ adjacent to it. The shaded pentagons are those across which the neighbour $B_j$ of $A$ is isometric to $A$. (b)~The same $2$-face $R$, regarded as a face of the neighbour $C$ of $A$ across $R$. The edge $Q_j\cap R$ of $P$ lies in exactly three facets, namely $A$, $B_j$ and $C$, so the face of $C$ across this edge is $B_j\cap C$, and by~(iii) it is a hexagon whenever $B_j\cong A$. Thus $R$ has at least $h$ hexagonal neighbours in $C$. The faces of $C$ across the other edges of $R$ (dotted, marked~?) are not determined.}\label{fig:f16b}
\end{figure}
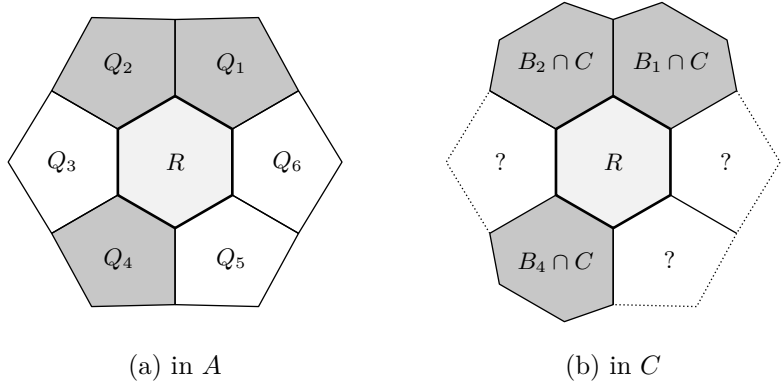

\begin{proof}
For $A\cong F_{15},F_{16a},F_{16b}$, a neighbour $B\cong A$ has $\dd(B,Q)\le2(m(A)-12)$ by~\eqref{eq:defect} and~\eqref{eq:euler}, and this is below the respective cutoff $13,10,10$. So all placements of $A$ on itself that survive Lemma~\ref{lem:placement}(i)--(iii) were found in the proof of Lemma~\ref{lem:single}. Two such placements, on faces $Q,Q'$ of $A$, are \emph{incompatible} if $Q=Q'$, or if $Q,Q'$ share an edge and at one of its endpoints $z$ the edges $t_B(z)$ and $t_C(z)$ of the two placements have different lengths. By Lemma~\ref{lem:placement}(iv), the placements realized by the neighbours of $A$ are pairwise compatible.

(i) The twelve faces of $F_{15}$ that admit a placement of $F_{15}$ fall into two sets of six, and within each set the incompatibilities contain a $6$-cycle (Appendix~\ref{subsec:simultaneous-data}). A $6$-cycle has no independent set with more than three vertices, so at most $3+3$ same-type neighbours occur.

(ii) For $F_{16a}$ the surviving placements fall into four classes of pairwise incompatible ones (Appendix~\ref{subsec:simultaneous-data}), so at most one placement from each class is realized.

(iii) One checks, for each of the $24$ incident pairs $(Q,R)$ and each surviving placement of $F_{16b}$ on $Q$, that the face of the target across the common edge $Q\cap R$ is a hexagon, and that the placement is incompatible with every surviving placement of $F_{16b}$ on $R$. For a realized placement, the face of $B$ across $Q\cap R$ is $B\cap C$.

(iv) If $h=0$, then $\dd(C,R)\ge7$ by Lemma~\ref{lem:single}(ii). Let $h\ge1$, and let $B_1,\dots,B_h$ be the neighbours isometric to $A$ across the pentagons $Q_1,\dots,Q_h$ adjacent to $R$. By~(iii), $C\not\cong A$, so $\dd(C,R)\ge10$ by Lemma~\ref{lem:single}(ii), and the faces $B_j\cap C$ are $h$ distinct hexagons of $C$ adjacent to $R$ (Figure~\ref{fig:f16b}). If $m(C)\le23$, an inspection of the catalogue shows the following: for every polyhedron $C'$ with at most $23$ faces and every hexagonal face $R'$ of $C'$ with $\dd(C',R')\ge10$ and at least $h$ hexagonal neighbours, one has $\dd(C',R')\ge10,10,10,10,11,14$ for $h=0,\dots,5$, respectively, and no such $R'$ exists for $h=6$. Each of these values is at least $2h+3$. If $m(C)\ge24$, then $C\not\cong L_6$ and Lemma~\ref{lem:collar} gives $\dd_{\ge2}(C,R)\ge2$, hence $\dd(C,R)\ge m(C)-10\ge14$ and if $h=6$, then $\dd_1(C,R)=6$, so~\eqref{eq:defect} gives $\dd_{\ge2}(C,R)=m(C)-19$ and $\dd(C,R)=2m(C)-31\ge17$. In all cases $\dd(C,R)\ge2h+3$.
\end{proof}

\subsection{L\"obell facets}\label{subsec:lobell}
We use the L\"obell polyhedra $L_n$ and their edge types $\ell_{n,1},\ell_{n,2},\ell_{n,3}$ defined in Section~\ref{subsec:unions} (Figure~\ref{fig:lobell}). For $n\ge6$ put $\tau_n=1+2\cos(2\pi/n)\ge2$. By Theorem~\ref{thm:andreev}, the rotational symmetry of $\Gamma(L_n)$ of order $n$ is induced by an isometry, so the bases are regular right-angled $n$-gons. Their sides satisfy $\cosh\ell_{n,1}=4\cos^2(\pi/n)-1=\tau_n$. Together with~\eqref{eq:pentagon} applied to a side face, this gives
\begin{equation}\label{eq:lobell}
 \cosh\ell_{n,1}=\tau_n,\qquad \cosh^2\ell_{n,2}=\frac{\tau_n}{\tau_n-1},\qquad \cosh^2\ell_{n,3}=\tau_n+1 .
\end{equation}
In particular the side faces are not regular, since $\tau_n\ne\varphi$. Every vertex of $L_n$ lies on exactly one edge of type $\ell_{n,2}$. For $L_6$ we have $\tau_6=2$, so $\cosh\ell_{6,1}=2$, $\cosh\ell_{6,2}=\sqrt2$ and $\cosh\ell_{6,3}=\sqrt3$.

\begin{lemma}\label{lem:lobell}
Let $A\cong L_n$ with $n\ge6$ and $B$ be neighbouring facets of a compact right-angled hyperbolic $4$-polytope $P$, and let $Q=A\cap B$. Then $\dd(B,Q)\ge4n-13$. Moreover, if $n=7$ and $B'$ is a neighbour of $A$ whose common face with $A$ contains an edge $e$ of type $\ell_{7,3}$, then the face of $B'$ through $e$ other than $A\cap B'$ is not a pentagon.
\end{lemma}

\begin{proof}
Let $e=zz'$ be an edge of $A\cap B$, and let $a,a'$ be the edges of type $\ell_{n,2}$ at $z$ and $z'$. By Lemma~\ref{lem:placement}(iii), if $\cosh\ell(e)\ge\cosh\ell(a)\cosh\ell(a')=\tau_n/(\tau_n-1)$, then the face of $B$ through $e$ other than $A\cap B$ is not a pentagon. For a base edge, $\cosh\ell(e)=\tau_n\ge \tau_n/(\tau_n-1)$ because $\tau_n\ge2$. So all $n$ faces of $B$ adjacent to $Q$ are nonpentagonal and $\dd_1(B,Q)\ge n$. In particular $B\not\cong L_n$, and Lemma~\ref{lem:collar} gives $\dd_{\ge2}(B,Q)\ge n-4$. By~\eqref{eq:defect},
\[
 \dd(B,Q)\ge(n-5)+n+2(n-4)=4n-13 .
\]
For $n=7$, the number $\tau=\tau_7$ is a root of $\tau^3-2\tau^2-\tau+1$, which is equivalent to $(\tau+1)(\tau-1)^2=\tau^2$. So an edge $e$ of type $\ell_{7,3}$ has $\cosh^2\ell(e)=\tau_7+1=\tau_7^2/(\tau_7-1)^2$, that is, $\cosh\ell(e)=\tau_7/(\tau_7-1)$, and Lemma~\ref{lem:placement}(iii), applied to $A$ and $B'$ with $a,a'$ the edges of type $\ell_{7,2}$ at the endpoints of $e$, gives the claim.
\end{proof}

\begin{corollary}\label{cor:L67}
If $N\le120$, then no facet of a compact right-angled hyperbolic $4$-polytope $P$ is isometric to a L\"obell polyhedron $L_6$ or $L_7$.
\end{corollary}

\begin{proof}
If $A\cong L_6$, Lemma~\ref{lem:single}(i) bounds each of the twelve pentagonal contributions to $\dd_P(A)$ from below by $10$, and Lemma~\ref{lem:lobell} bounds each of the two base contributions from below by $11$. By Proposition~\ref{prop:identity},
\[
 2N\ge30+5\cdot14+12\cdot10+2\cdot11=242 ,
\]
so $N\ge121$. For $A\cong L_7$ we get in the same way $2N\ge30+5\cdot16+14\cdot8+2\cdot15=252$, so $N\ge126$.
\end{proof}

The second assertion of Lemma~\ref{lem:lobell}, for neighbours across side pentagons, is one of the conditions used in the finite check of Lemma~\ref{lem:single}(i) for $L_7$.

\subsection{Proof of Proposition~\ref{prop:small}}\label{subsec:small-proof}
\begin{proof}
Assume $N\le120$. Corollary~\ref{cor:L67} excludes $L_6$ and $L_7$. For $F_{15}$, Lemma~\ref{lem:simultaneous}(i) and Lemma~\ref{lem:single}(ii) give at most six neighbours of defect at least $4$ and at least nine of defect at least $13$, so $\dd_P(A)\ge6\cdot4+9\cdot13=141$, and~\eqref{eq:identity} gives $N\ge123$. The cases $F_{16a}$ and $F_{17a},F_{17b},F_{17c}$ are entirely similar, using Lemma~\ref{lem:single} and, for $F_{16a}$, Lemma~\ref{lem:simultaneous}(ii), and the resulting bounds are collected in Table~\ref{tab:five}.

\begin{table}[ht]
\centering
\caption{Lower bounds for $\dd_P(A)$ and $N$.}\label{tab:five}
\begin{tabular}{lcc}\toprule
$A$ & lower bound for $\dd_P(A)$ & lower bound for $N$\\\midrule
$F_{15}$  &$6\cdot4+9\cdot13=141$&$123$\\
$F_{16a}$ &$4\cdot6+12\cdot10=144$&$127$\\
$F_{17a}$ &$6\cdot7+11\cdot8=130$&$123$\\
$F_{17b}$ &$3\cdot7+14\cdot8=133$&$124$\\
$F_{17c}$ &$2\cdot7+15\cdot8=134$&$125$\\\bottomrule
\end{tabular}
\end{table}

For $F_{16b}$, a neighbour of small defect across a pentagon is compensated by large defects across the adjacent hexagons. Let $A\cong F_{16b}$, let $s$ be the number of pentagons of $A$ across which the neighbour is isometric to $A$, and for the four hexagons $R_1,\dots,R_4$ of $A$ let $h_i$ be the number of such pentagons adjacent to $R_i$ and $C_i$ the neighbour across $R_i$. Since every pentagon of $A$ meets exactly two hexagons, $\sum_ih_i=2s$. By Lemma~\ref{lem:single}(ii) the pentagonal contributions to $\dd_P(A)$ add up to at least $6s+10(12-s)$, and by Lemma~\ref{lem:simultaneous}(iv)
\[
 \dd_P(A)\ge120-4s+\sum_{i=1}^4(2h_i+3)=132 .
\]
Hence $2N\ge30+80+132$ and $N\ge121$.

This excludes every type with at most $17$ faces except $D$ and $E$. Consequently every neighbour $B\not\cong D$ of a facet of $P$ has at least $17$ faces, and its defect at any face $Q$ is at least $7$: for $m(B)=17$ this is Lemma~\ref{lem:small-facts}(ii), and for $m(B)\ge18$ it follows from $\dd(B,Q)\ge m(B)-11$.

Let $A$ be a facet with $18$ or $19$ faces. By Lemma~\ref{lem:small-facts}(i), the facet $A$ has no regular pentagon, so no neighbour of $A$ is a dodecahedron, and all $m(A)$ neighbours have defect at least $7$. Then~\eqref{eq:identity} gives $2N\ge30+12\,m(A)$, so $N\ge15+6\cdot18=123$. The same argument applies to $L_n$ with $n\ge8$, which has $2n+2\ge18$ faces and, by~\eqref{eq:lobell}, no regular pentagon. This completes the proof of Proposition~\ref{prop:small}.
\end{proof}

The case $A\cong F_{16b}$ is the tightest one in the proof of Proposition~\ref{prop:small}. It is also the only case in which we do not bound the number of neighbours isometric to $A$. By Table~\ref{tab:same-type}, placements of $F_{16b}$ on itself survive on all twelve pentagons and on all four hexagons. What makes the argument work is that a cheap neighbour across a pentagon makes the neighbours across both adjacent hexagons more expensive. Since each pentagon meets exactly two hexagons, the saving of $4$ at a pentagon is exactly offset by the two increases of $2$ given by Lemma~\ref{lem:simultaneous}(iv), and the bound $\dd_P(A)\ge132$ does not depend on $s$.

\section{Regular pentagons and a bound on the size of facets}\label{sec:bound}

Throughout this section we assume $A$ is a facet polyhedron of a compact right-angled hyperbolic $4$-polytope $P$ with $N\le120$. By Proposition~\ref{prop:small} and~\eqref{eq:defect}, every neighbour of $A$ that is not a right-angled dodecahedron has defect at least $7$. Dodecahedral neighbours of $A$, on the other hand, contribute nothing to $\dd_P(A)$, and they can only be attached along regular pentagons of $A$. We show in this section that a facet cannot have too many regular pentagons.

\begin{lemma}\label{lem:firstcost}
Suppose $N\le120$. Let $B$ be a facet of $P$ that is not isometric to the right-angled dodecahedron $D$, and let $Q$ be a face of $B$. Then
\[
 \dd(B,Q)\ge
 \begin{cases}
  8&\text{if $Q$ is a regular pentagon},\\
  9&\text{if $Q$ is a nonregular pentagon},\\
  2k(Q)-5&\text{if $k(Q)\ge6$}.
 \end{cases}
\]
\end{lemma}

\begin{proof}
If $B$ is isometric to the double $E$ of the right-angled dodecahedron, all pentagons of $B$ are regular, and the defects are $8$ or $10$ at pentagons and $7=2\cdot6-5$ at hexagons. Otherwise, by Proposition~\ref{prop:small}, $m(B)\ge20$ and $B$ is not a L\"obell polyhedron, so Lemma~\ref{lem:collar} gives $\dd(B,Q)=m(B)-12+\dd_{\ge2}(B,Q)\ge8+(k-4)=k+4$, where $k=k(Q)$. This is at least $9$ for $k=5$ and at least $2k-5$ for $k=6,7$. For $k\ge8$, \eqref{eq:collar-bound} gives $\dd(B,Q)\ge(2k-9)+(k-4)=3k-13\ge2k-5$.
\end{proof}

For a polyhedron $B$, let $U(B)$ be the set of its regular pentagonal faces, regarded as a set of vertices of $\Gamma(B)$, and let $\Gamma(B)[U(B)]$ be the subgraph it induces.

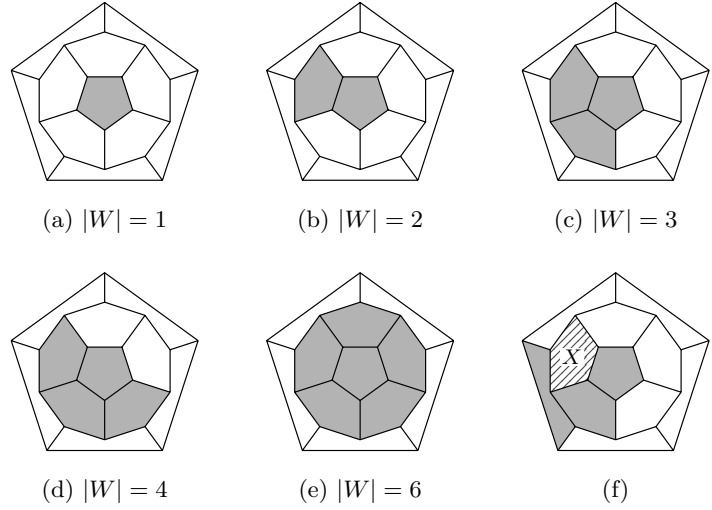
\begin{figure}[ht]
\centering
\begin{tikzpicture}[scale=1.3,line width=0.45pt,line join=round]
\begin{scope}[xshift=0.00cm,yshift=-0.00cm]
  \fill[black!30] (-0.176,0.243) -- (-0.285,-0.093) -- (-0.000,-0.300) -- (0.285,-0.093) -- (0.176,0.243) -- cycle;
  \draw (0.000,1.000) -- (-0.951,0.309) (0.000,1.000) -- (0.000,0.700) (0.000,0.700) -- (-0.411,0.566) (-0.411,0.566) -- (-0.666,0.216) (-0.411,0.566) -- (-0.176,0.243) (-0.176,0.243) -- (-0.285,-0.093) (-0.951,0.309) -- (-0.588,-0.809) (-0.951,0.309) -- (-0.666,0.216) (-0.666,0.216) -- (-0.666,-0.216) (-0.666,-0.216) -- (-0.411,-0.566) (-0.666,-0.216) -- (-0.285,-0.093) (-0.285,-0.093) -- (-0.000,-0.300) (-0.588,-0.809) -- (0.588,-0.809) (-0.588,-0.809) -- (-0.411,-0.566) (-0.411,-0.566) -- (-0.000,-0.700) (-0.000,-0.700) -- (0.411,-0.566) (-0.000,-0.700) -- (-0.000,-0.300) (-0.000,-0.300) -- (0.285,-0.093) (0.588,-0.809) -- (0.951,0.309) (0.588,-0.809) -- (0.411,-0.566) (0.411,-0.566) -- (0.666,-0.216) (0.666,-0.216) -- (0.666,0.216) (0.666,-0.216) -- (0.285,-0.093) (0.285,-0.093) -- (0.176,0.243) (0.951,0.309) -- (0.000,1.000) (0.951,0.309) -- (0.666,0.216) (0.666,0.216) -- (0.411,0.566) (0.411,0.566) -- (0.000,0.700) (0.411,0.566) -- (0.176,0.243) (0.176,0.243) -- (-0.176,0.243);
  \node[font=\footnotesize] at (0,-1.22) {(a) $|W|=1$};
\end{scope}
\begin{scope}[xshift=2.60cm,yshift=-0.00cm]
  \fill[black!30] (-0.176,0.243) -- (-0.285,-0.093) -- (-0.000,-0.300) -- (0.285,-0.093) -- (0.176,0.243) -- cycle;
  \fill[black!30] (-0.411,0.566) -- (-0.176,0.243) -- (-0.285,-0.093) -- (-0.666,-0.216) -- (-0.666,0.216) -- cycle;
  \draw (0.000,1.000) -- (-0.951,0.309) (0.000,1.000) -- (0.000,0.700) (0.000,0.700) -- (-0.411,0.566) (-0.411,0.566) -- (-0.666,0.216) (-0.411,0.566) -- (-0.176,0.243) (-0.176,0.243) -- (-0.285,-0.093) (-0.951,0.309) -- (-0.588,-0.809) (-0.951,0.309) -- (-0.666,0.216) (-0.666,0.216) -- (-0.666,-0.216) (-0.666,-0.216) -- (-0.411,-0.566) (-0.666,-0.216) -- (-0.285,-0.093) (-0.285,-0.093) -- (-0.000,-0.300) (-0.588,-0.809) -- (0.588,-0.809) (-0.588,-0.809) -- (-0.411,-0.566) (-0.411,-0.566) -- (-0.000,-0.700) (-0.000,-0.700) -- (0.411,-0.566) (-0.000,-0.700) -- (-0.000,-0.300) (-0.000,-0.300) -- (0.285,-0.093) (0.588,-0.809) -- (0.951,0.309) (0.588,-0.809) -- (0.411,-0.566) (0.411,-0.566) -- (0.666,-0.216) (0.666,-0.216) -- (0.666,0.216) (0.666,-0.216) -- (0.285,-0.093) (0.285,-0.093) -- (0.176,0.243) (0.951,0.309) -- (0.000,1.000) (0.951,0.309) -- (0.666,0.216) (0.666,0.216) -- (0.411,0.566) (0.411,0.566) -- (0.000,0.700) (0.411,0.566) -- (0.176,0.243) (0.176,0.243) -- (-0.176,0.243);
  \node[font=\footnotesize] at (0,-1.22) {(b) $|W|=2$};
\end{scope}
\begin{scope}[xshift=5.20cm,yshift=-0.00cm]
  \fill[black!30] (-0.176,0.243) -- (-0.285,-0.093) -- (-0.000,-0.300) -- (0.285,-0.093) -- (0.176,0.243) -- cycle;
  \fill[black!30] (-0.411,0.566) -- (-0.176,0.243) -- (-0.285,-0.093) -- (-0.666,-0.216) -- (-0.666,0.216) -- cycle;
  \fill[black!30] (-0.666,-0.216) -- (-0.285,-0.093) -- (-0.000,-0.300) -- (-0.000,-0.700) -- (-0.411,-0.566) -- cycle;
  \draw (0.000,1.000) -- (-0.951,0.309) (0.000,1.000) -- (0.000,0.700) (0.000,0.700) -- (-0.411,0.566) (-0.411,0.566) -- (-0.666,0.216) (-0.411,0.566) -- (-0.176,0.243) (-0.176,0.243) -- (-0.285,-0.093) (-0.951,0.309) -- (-0.588,-0.809) (-0.951,0.309) -- (-0.666,0.216) (-0.666,0.216) -- (-0.666,-0.216) (-0.666,-0.216) -- (-0.411,-0.566) (-0.666,-0.216) -- (-0.285,-0.093) (-0.285,-0.093) -- (-0.000,-0.300) (-0.588,-0.809) -- (0.588,-0.809) (-0.588,-0.809) -- (-0.411,-0.566) (-0.411,-0.566) -- (-0.000,-0.700) (-0.000,-0.700) -- (0.411,-0.566) (-0.000,-0.700) -- (-0.000,-0.300) (-0.000,-0.300) -- (0.285,-0.093) (0.588,-0.809) -- (0.951,0.309) (0.588,-0.809) -- (0.411,-0.566) (0.411,-0.566) -- (0.666,-0.216) (0.666,-0.216) -- (0.666,0.216) (0.666,-0.216) -- (0.285,-0.093) (0.285,-0.093) -- (0.176,0.243) (0.951,0.309) -- (0.000,1.000) (0.951,0.309) -- (0.666,0.216) (0.666,0.216) -- (0.411,0.566) (0.411,0.566) -- (0.000,0.700) (0.411,0.566) -- (0.176,0.243) (0.176,0.243) -- (-0.176,0.243);
  \node[font=\footnotesize] at (0,-1.22) {(c) $|W|=3$};
\end{scope}
\begin{scope}[xshift=0.00cm,yshift=-2.75cm]
  \fill[black!30] (-0.176,0.243) -- (-0.285,-0.093) -- (-0.000,-0.300) -- (0.285,-0.093) -- (0.176,0.243) -- cycle;
  \fill[black!30] (-0.411,0.566) -- (-0.176,0.243) -- (-0.285,-0.093) -- (-0.666,-0.216) -- (-0.666,0.216) -- cycle;
  \fill[black!30] (-0.666,-0.216) -- (-0.285,-0.093) -- (-0.000,-0.300) -- (-0.000,-0.700) -- (-0.411,-0.566) -- cycle;
  \fill[black!30] (-0.000,-0.700) -- (-0.000,-0.300) -- (0.285,-0.093) -- (0.666,-0.216) -- (0.411,-0.566) -- cycle;
  \draw (0.000,1.000) -- (-0.951,0.309) (0.000,1.000) -- (0.000,0.700) (0.000,0.700) -- (-0.411,0.566) (-0.411,0.566) -- (-0.666,0.216) (-0.411,0.566) -- (-0.176,0.243) (-0.176,0.243) -- (-0.285,-0.093) (-0.951,0.309) -- (-0.588,-0.809) (-0.951,0.309) -- (-0.666,0.216) (-0.666,0.216) -- (-0.666,-0.216) (-0.666,-0.216) -- (-0.411,-0.566) (-0.666,-0.216) -- (-0.285,-0.093) (-0.285,-0.093) -- (-0.000,-0.300) (-0.588,-0.809) -- (0.588,-0.809) (-0.588,-0.809) -- (-0.411,-0.566) (-0.411,-0.566) -- (-0.000,-0.700) (-0.000,-0.700) -- (0.411,-0.566) (-0.000,-0.700) -- (-0.000,-0.300) (-0.000,-0.300) -- (0.285,-0.093) (0.588,-0.809) -- (0.951,0.309) (0.588,-0.809) -- (0.411,-0.566) (0.411,-0.566) -- (0.666,-0.216) (0.666,-0.216) -- (0.666,0.216) (0.666,-0.216) -- (0.285,-0.093) (0.285,-0.093) -- (0.176,0.243) (0.951,0.309) -- (0.000,1.000) (0.951,0.309) -- (0.666,0.216) (0.666,0.216) -- (0.411,0.566) (0.411,0.566) -- (0.000,0.700) (0.411,0.566) -- (0.176,0.243) (0.176,0.243) -- (-0.176,0.243);
  \node[font=\footnotesize] at (0,-1.22) {(d) $|W|=4$};
\end{scope}
\begin{scope}[xshift=2.60cm,yshift=-2.75cm]
  \fill[black!30] (-0.176,0.243) -- (-0.285,-0.093) -- (-0.000,-0.300) -- (0.285,-0.093) -- (0.176,0.243) -- cycle;
  \fill[black!30] (-0.411,0.566) -- (-0.176,0.243) -- (-0.285,-0.093) -- (-0.666,-0.216) -- (-0.666,0.216) -- cycle;
  \fill[black!30] (-0.666,-0.216) -- (-0.285,-0.093) -- (-0.000,-0.300) -- (-0.000,-0.700) -- (-0.411,-0.566) -- cycle;
  \fill[black!30] (-0.000,-0.700) -- (-0.000,-0.300) -- (0.285,-0.093) -- (0.666,-0.216) -- (0.411,-0.566) -- cycle;
  \fill[black!30] (0.666,-0.216) -- (0.285,-0.093) -- (0.176,0.243) -- (0.411,0.566) -- (0.666,0.216) -- cycle;
  \fill[black!30] (0.411,0.566) -- (0.176,0.243) -- (-0.176,0.243) -- (-0.411,0.566) -- (0.000,0.700) -- cycle;
  \draw (0.000,1.000) -- (-0.951,0.309) (0.000,1.000) -- (0.000,0.700) (0.000,0.700) -- (-0.411,0.566) (-0.411,0.566) -- (-0.666,0.216) (-0.411,0.566) -- (-0.176,0.243) (-0.176,0.243) -- (-0.285,-0.093) (-0.951,0.309) -- (-0.588,-0.809) (-0.951,0.309) -- (-0.666,0.216) (-0.666,0.216) -- (-0.666,-0.216) (-0.666,-0.216) -- (-0.411,-0.566) (-0.666,-0.216) -- (-0.285,-0.093) (-0.285,-0.093) -- (-0.000,-0.300) (-0.588,-0.809) -- (0.588,-0.809) (-0.588,-0.809) -- (-0.411,-0.566) (-0.411,-0.566) -- (-0.000,-0.700) (-0.000,-0.700) -- (0.411,-0.566) (-0.000,-0.700) -- (-0.000,-0.300) (-0.000,-0.300) -- (0.285,-0.093) (0.588,-0.809) -- (0.951,0.309) (0.588,-0.809) -- (0.411,-0.566) (0.411,-0.566) -- (0.666,-0.216) (0.666,-0.216) -- (0.666,0.216) (0.666,-0.216) -- (0.285,-0.093) (0.285,-0.093) -- (0.176,0.243) (0.951,0.309) -- (0.000,1.000) (0.951,0.309) -- (0.666,0.216) (0.666,0.216) -- (0.411,0.566) (0.411,0.566) -- (0.000,0.700) (0.411,0.566) -- (0.176,0.243) (0.176,0.243) -- (-0.176,0.243);
  \node[font=\footnotesize] at (0,-1.22) {(e) $|W|=6$};
\end{scope}
\begin{scope}[xshift=5.20cm,yshift=-2.75cm]
  \fill[black!30] (-0.176,0.243) -- (-0.285,-0.093) -- (-0.000,-0.300) -- (0.285,-0.093) -- (0.176,0.243) -- cycle;
  \fill[black!30] (-0.666,-0.216) -- (-0.285,-0.093) -- (-0.000,-0.300) -- (-0.000,-0.700) -- (-0.411,-0.566) -- cycle;
  \fill[black!30] (-0.951,0.309) -- (-0.588,-0.809) -- (-0.411,-0.566) -- (-0.666,-0.216) -- (-0.666,0.216) -- cycle;
  \fill[pattern=north east lines,pattern color=black!70] (-0.411,0.566) -- (-0.176,0.243) -- (-0.285,-0.093) -- (-0.666,-0.216) -- (-0.666,0.216) -- cycle;
  \draw (0.000,1.000) -- (-0.951,0.309) (0.000,1.000) -- (0.000,0.700) (0.000,0.700) -- (-0.411,0.566) (-0.411,0.566) -- (-0.666,0.216) (-0.411,0.566) -- (-0.176,0.243) (-0.176,0.243) -- (-0.285,-0.093) (-0.951,0.309) -- (-0.588,-0.809) (-0.951,0.309) -- (-0.666,0.216) (-0.666,0.216) -- (-0.666,-0.216) (-0.666,-0.216) -- (-0.411,-0.566) (-0.666,-0.216) -- (-0.285,-0.093) (-0.285,-0.093) -- (-0.000,-0.300) (-0.588,-0.809) -- (0.588,-0.809) (-0.588,-0.809) -- (-0.411,-0.566) (-0.411,-0.566) -- (-0.000,-0.700) (-0.000,-0.700) -- (0.411,-0.566) (-0.000,-0.700) -- (-0.000,-0.300) (-0.000,-0.300) -- (0.285,-0.093) (0.588,-0.809) -- (0.951,0.309) (0.588,-0.809) -- (0.411,-0.566) (0.411,-0.566) -- (0.666,-0.216) (0.666,-0.216) -- (0.666,0.216) (0.666,-0.216) -- (0.285,-0.093) (0.285,-0.093) -- (0.176,0.243) (0.951,0.309) -- (0.000,1.000) (0.951,0.309) -- (0.666,0.216) (0.666,0.216) -- (0.411,0.566) (0.411,0.566) -- (0.000,0.700) (0.411,0.566) -- (0.176,0.243) (0.176,0.243) -- (-0.176,0.243);
  \node[fill=white,inner sep=0.6pt,font=\scriptsize] at (-0.445,0.144) {$X$};
  \node[font=\footnotesize] at (0,-1.22) {(f)};
\end{scope}
\end{tikzpicture}
\caption{Connected sets $W$ of faces of the dodecahedron (shaded) such that no other face has three consecutive neighbours in $W$ (proof of Lemma~\ref{lem:clusters}). Up to symmetry there are exactly five: (a)~one face, (b)~two adjacent faces, (c)~the three faces at a vertex, (d)~the four faces at an edge, and (e)~a face together with its five neighbours. (f)~An excluded set: the hatched face $X$ has three consecutive neighbours in $W$, so in $B$ it would have three consecutive sides of length $\ell_0$.}\label{fig:clusters}
\end{figure}

\begin{lemma}\label{lem:clusters}
Let $B\not\cong D$ be a compact right-angled polyhedron. Every connected component of $\Gamma(B)[U(B)]$ has $1$, $2$, $3$, $4$ or $6$ vertices and is adjacent to at least five faces of $B$ outside $U(B)$. If $B$ has $m$ faces, $r=|U(B)|$ regular pentagons, and the graph $\Gamma(B)[U(B)]$ has $s$ components, then
\begin{equation}\label{eq:clusters}
 5r\le4m-8.
\end{equation}
\end{lemma}
\begin{proof}
For a regular pentagonal face $Q$ of $B$, let $D_Q$ be the copy of $D$ that has $Q$ as a face and lies on the same side of the plane of $Q$ as $B$. Let $Q,Q'\in U(B)$ share an edge $e$. Both $Q'$ and the face of $D_Q$ across $e$ are regular right-angled pentagons with side $e$, lying in the plane through $e$ orthogonal to $Q$ and on the same side of $Q$ as $B$, and hence they coincide. So $Q'$ is a face of $D_Q$, and $D_Q$ lies on the same side of the plane of $Q'$ as $B$, since both contain $Q$. Thus $D_{Q'}=D_Q$, and $D_Q$ is the same copy $D_0$ for all faces $Q$ of a component $W$. In particular, the component $W$ is a connected set of faces of $D_0$.

Let $X$ be a face of $D_0$ not in $W$ but adjacent to $W$. Since all dihedral angles are right, the face $Y$ of $B$ across the corresponding edge lies in the plane of $X$. So distinct $X$ give distinct $Y$, and $Y\notin U(B)$ because $W$ is a component. If $X$ had three consecutive neighbours in $W$, then $Y$ would have three consecutive sides of length $\ell_0$ and would be a regular pentagon by Lemma~\ref{lem:pentagon}, a contradiction (Figure~\ref{fig:clusters}(f)). Checking all $2^{12}-1$ nonempty sets of faces of $D_0$, one finds the following connected sets with this property:
\begin{center}
\begin{tabular}{lcccccc}\toprule
$|W|$ &1&2&3&4&6&12\\
number of such sets&12&30&20&30&12&1\\
faces adjacent to $W$&5&6&6&6&5&0\\\bottomrule
\end{tabular}
\end{center}
Up to symmetry these are the sets of Figure~\ref{fig:clusters}(a)--(e) and the set of all twelve faces. The last is impossible: it would give $\partial D_0\subseteq\partial B$, hence $B=D_0$.

Finally, let $s$ be the number of components and $b=m-r$. Contract each component in the planar graph $\Gamma(B)$ and keep only the edges between components and faces outside $U(B)$. The result is a simple planar bipartite graph with $s+b$ vertices and at least $5s$ edges, so $5s\le2(s+b)-4$, that is, $3s\le2b-4$. Since $r\le6s$, we get $r\le4b-8$, that is, $5r\le4m-8$.
\end{proof}
\begin{proposition}\label{prop:upper}
If $N\le120$, then every facet of $P$ is isometric to $D$ or $E$, or has between $20$ and $27$ faces.
\end{proposition}

\begin{proof}
Let $A\not\cong D$ be a facet with $m$ faces, $p$ pentagons and $r$ regular pentagons, and let $d$ be the number of its dodecahedral neighbours. Then $d\le r\le p$, since a dodecahedral neighbour can only be attached along a regular pentagon. Apply Lemma~\ref{lem:firstcost} to all neighbours, using the bound $8=(2\cdot5-5)+3$ for every pentagonal face and then subtracting $8$ for each dodecahedral neighbour. Since $\sum_Qk(Q)=6m-12$,
\[
 \dd_P(A)\ge\sum_{Q}(2k(Q)-5)+3p-8d=7m-24+3p-8d .
\]
With~\eqref{eq:identity} and~\eqref{eq:clusters} this gives
\[
 N\ge3+6m+\tfrac32p-4d\ge3+6m-\tfrac52r\ge7+4m .
\]
Hence $m\le28$. If $m=28$, then $r\le\lfloor(4m-8)/5\rfloor=20$, and the middle inequality gives $N\ge3+168-50=121$. So $m\le27$, and Proposition~\ref{prop:small} completes the proof.
\end{proof}

For later use we record a slightly sharper form of the estimate in the proof of Proposition~\ref{prop:upper}, which keeps the extra unit that Lemma~\ref{lem:firstcost} provides at each of the $p-r$ nonregular pentagons.

\begin{lemma}\label{lem:R1}
Suppose $N\le120$, and let $A\not\cong D$ be a facet of $P$ with $m$ faces, $p$ pentagons, $r$ regular pentagons and $d$ dodecahedral neighbours. Then
\begin{equation}\label{eq:R1}
 240\ge2N\ge6+12m+4p-r-8d\ge6+12m+4p-9r .
\end{equation}
\end{lemma}

The proof is the same as that of Proposition~\ref{prop:upper}, with the bound $\dd_P(A)\ge7m-24+4p-r-8d$ in place of $\dd_P(A)\ge7m-24+3p-8d$, and we omit it. Thus a facet with many faces must have many regular pentagons. In the next section we show that this cannot happen except for $\Tt$ and $\Th$.

\section{The facets of a polytope with at most 120 facets}\label{sec:classification}

We continue to assume $N\le120$. By Proposition~\ref{prop:upper}, it remains to determine the facets with $20$ to $27$ faces.

\begin{proposition}\label{prop:classification}
Every facet of a compact right-angled polytope $P\subset\HH^4$ with at most $120$ facets is isometric to $D$, $E$, $\Tt$ or $\Th$.
\end{proposition}

\subsection{Candidate sets of regular pentagons}\label{subsec:candidates}
Inequality~\eqref{eq:R1} only involves the combinatorial type of a facet and its set of regular pentagons. We therefore enumerate, for each combinatorial type, the sets of faces that could be its set of regular pentagons, and then test them metrically. The following necessary conditions make the enumeration finite and small.

\begin{lemma}\label{lem:regular-sets}
Let $B\not\cong D$ be a compact right-angled polyhedron, $G=\Gamma(B)$ and $U=U(B)$. Then:
\begin{enumerate}
\item every vertex of $U$ has degree five in $G$;
\item every vertex of degree five not in $U$ has at most one neighbour in $U$;
\item no vertex of degree at least six has three consecutive neighbours (in the cyclic order of its link) in $U$;
\item every connected component of $G[U]$ has $1$, $2$, $3$, $4$ or $6$ vertices;
\item $U$ is invariant under all automorphisms of $G$.
\end{enumerate}
\end{lemma}

\begin{proof}
Part~(i) is clear. A pentagon adjacent to two regular pentagons has two sides of length $\ell_0$ and is regular by Lemma~\ref{lem:pentagon}(i), which gives~(ii). Part~(iii) is Lemma~\ref{lem:pentagon}(ii), (iv) is part of Lemma~\ref{lem:clusters}, and~(v) follows from Theorem~\ref{thm:andreev}.
\end{proof}

It is important here that $U(B)$ is the set of \emph{all} regular pentagons: a proper subset need not satisfy~(ii) or~(v).

For each $m$ let $\mathcal G_m$ be the catalogue of $5$-connected triangulations of the sphere with $m$ vertices (Appendix~\ref{subsec:catalogues}). By Lemma~\ref{lem:5conn} it contains $\Gamma(A)$ for every facet $A$ with $m$ faces, up to isomorphism. Given $G\in\mathcal G_m$ and an integer $R$, a \emph{candidate set} of size at least $R$ is a set $U\subseteq V(G)$ with $|U|\ge R$ satisfying conditions~(i)--(v) of Lemma~\ref{lem:regular-sets}. If $A$ is a facet of $P$ with $\Gamma(A)\cong G$ and if~\eqref{eq:R1} or~\eqref{eq:R2} below forces $|U(A)|\ge R$, then $U(A)$ corresponds to a candidate set. To exclude $G$ it therefore suffices to show, for each candidate set $U$, that some face in $U$ is not a regular pentagon in the right-angled realization of $G$, that is, that one of its sides has length different from $\ell_0$. This is checked with the rigorous enclosures of Appendix~\ref{app:computation}. A worked example is given in Appendix~\ref{subsec:example}.

\subsection{Facets with twenty to twenty-three faces}
Let $p$ be the number of vertices of degree five in $G\in\mathcal G_m$. By~\eqref{eq:R1}, a facet $A$ with $\Gamma(A)\cong G$ has
\[
 |U(A)|\ge R_1(G):=\ceil{\frac{6+12m+4p-240}{9}}.
\]

\begin{lemma}\label{lem:stage1}
Every facet of $P$ with $20$ to $23$ faces is isometric to $\Tt$ or $\Th$.
\end{lemma}

\begin{table}[ht]
\centering
\caption{The search for Lemma~\ref{lem:stage1}: candidate sets of size at least $R_1(G)$.}\label{tab:stage1}
\begin{tabular}{r>{\raggedleft\arraybackslash}p{1.25cm}ccc}\toprule
$m$ & $|\mathcal G_m|$ & graphs with a candidate set & candidate sets & graphs not excluded\\\midrule
20&71&15&36&0\\
21&187&15&34&0\\
22&627&25&40&2\\
23&1970&23&28&0\\\bottomrule
\end{tabular}
\end{table}

\begin{proof}
Let $A$ be a facet of $P$ with $20\le m\le23$ faces, and let $G\in\mathcal G_m$ be isomorphic to $\Gamma(A)$. By Lemma~\ref{lem:regular-sets} and the bound above, the set $U(A)$ corresponds to a candidate set of size at least $R_1(G)$ in $G$. An exhaustive search shows that $78$ graphs in $\mathcal G_{20},\dots,\mathcal G_{23}$ carry such a candidate set (Table~\ref{tab:stage1}). For each of them we have rigorous enclosures of the normals of its right-angled realization, and for $76$ of them every candidate set contains a face with a side $e$ for which the enclosure of $\cosh\ell(e)$ does not contain $\varphi$. So $G$ is one of the two remaining graphs, the catalogue entries $(22,622)$ and $(22,626)$. Starting from $\Gamma(D)$, one constructs $\Gamma(E)$ by gluing along a face $Q_1$, and then $\Gamma(\Tt)$ and $\Gamma(\Th)$ by gluing along a face $Q_2$ at distance $2$, respectively $3$, from $Q_1$ (Section~\ref{subsec:unions}). These graphs are isomorphic to entries $(22,622)$ and $(22,626)$, respectively. Since $\Tt$ and $\Th$ are right-angled, Theorem~\ref{thm:andreev} shows that $A$ is isometric to $\Tt$ or $\Th$.
\end{proof}

We have shown the following.

\begin{corollary}\label{cor:intermediate}
If $N\le120$, then every facet of $P$ belongs to $\cA$ or has between $24$ and $27$ faces.
\end{corollary}

\subsection{Facets with twenty-four to twenty-seven faces}
Corollary~\ref{cor:intermediate} improves the lower bounds for defects of neighbours, and with them the required number of regular pentagons. Let $A$ be a facet with $24\le m(A)\le27$ and $B$ the neighbour across a face $Q$ of $A$, with $k=k(Q)$.

If $Q$ is a nonregular pentagon, then $B\notin\cA$, since all pentagons of members of $\cA$ are regular. So $m(B)\ge24$ by Corollary~\ref{cor:intermediate}, and $\dd(B,Q)\ge12+1=13$ by Lemma~\ref{lem:collar}. If $k\ge6$, we claim $\dd(B,Q)\ge g(k)$, where
\begin{equation}\label{eq:g}
 g(6)=7,\qquad g(7)=14,\qquad g(k)=\max\{k+8,\,3k-13\}\quad(k\ge8).
\end{equation}
Indeed, the neighbour $B$ is either one of $E,\Tt,\Th$ or has at least $24$ faces. The defects of $E$ at its hexagons are $7$. The defects of $\Tt$ are $15$ or $17$ at its hexagons and $14$ at its heptagons, and those of $\Th$ are $15$ at its hexagons. If $m(B)\ge24$, then $B$ is not a L\"obell polyhedron (Proposition~\ref{prop:small}) and Lemma~\ref{lem:collar} gives $\dd(B,Q)\ge12+(k-4)$ as well as $\dd(B,Q)\ge(2k-9)+(k-4)$. Since no member of $\cA$ has a face with more than seven sides, this proves the claim. Note that this estimate remains valid if other facets of $P$ have $24$ to $27$ faces.

Discarding the contributions of the regular pentagons, \eqref{eq:identity} now gives
\begin{equation}\label{eq:R2}
 |U(A)|\ge R_2(G):=\ceil{\frac{30+5m+13p+\sum_{k(Q)\ge6}g(k(Q))-240}{13}},
\end{equation}
where $G\cong\Gamma(A)$, $m=m(A)$, $p$ is the number of pentagons of $A$, and the sum runs over the faces of $A$ with at least six sides. Like $R_1(G)$, this number depends only on the graph $G$. Since $A$ is not a L\"obell polyhedron, Lemma~\ref{lem:collar} gives $k(Q)\le\lfloor(m-3)/2\rfloor\le12$ for every face $Q$ of $A$. Accordingly, the L\"obell graphs $\Gamma(L_{11})\in\mathcal G_{24}$ and $\Gamma(L_{12})\in\mathcal G_{26}$ are left out of the search. All other graphs in $\mathcal G_m$, $24\le m\le27$, have maximal degree at most $\lfloor(m-3)/2\rfloor$, as they must by Lemma~\ref{lem:collar}.

\begin{figure}[ht]
\centering
\input{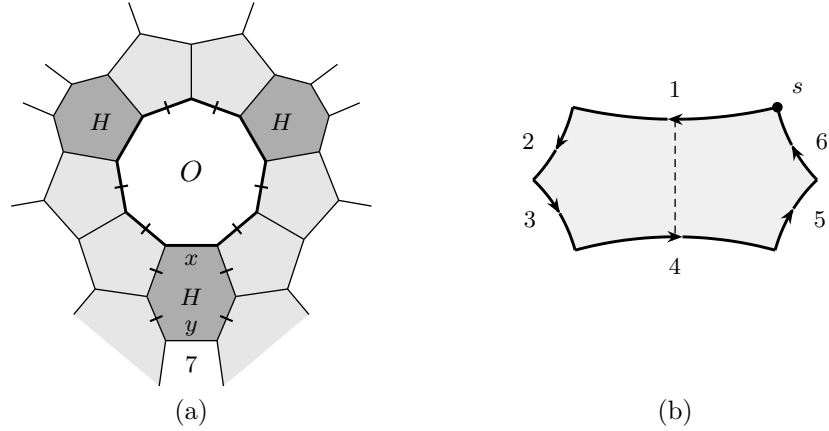}
\caption{The exclusion of $(25,23383)$ in Lemma~\ref{lem:stage2}. (a)~The faces around the nonagon $O$: three hexagons $H$ (dark) and six regular pentagons (light). The hexagon at the bottom is shown with its three further neighbours, two regular pentagons and a heptagon (marked $7$). Ticks mark the sides of $O$ and of this hexagon that are shared with regular pentagons, and these sides have length $\ell_0$. So the sides of this hexagon have lengths $(x,\ell_0,\ell_0,y,\ell_0,\ell_0)$, and~\eqref{eq:hexagon} gives $x=y=2\ell_0$. (b)~A right-angled path with side lengths $2\ell_0,\ell_0,\ell_0,2\ell_0,\ell_0,\ell_0$ (sides numbered in order, starting at $s$), drawn in the Poincar\'e disk, returns to its starting point $s$ after six sides. It bounds two regular pentagons glued along the dashed segment. The boundary of $O$ would have to follow this path, which is impossible for a nonagon.}\label{fig:nonagon}
\end{figure}

\begin{lemma}\label{lem:stage2}
No facet of $P$ has $24$ to $27$ faces.
\end{lemma}

\begin{table}[ht]
\centering
\caption{The search for Lemma~\ref{lem:stage2}: candidate sets of size at least $R_2(G)$.}\label{tab:stage2}
\begin{tabular}{r>{\raggedleft\arraybackslash}p{1.25cm}ccc}\toprule
$m$ & $|\mathcal G_m|$ & graphs with a candidate set & candidate sets & graphs not excluded\\\midrule
24&6833&10&10&0\\
25&23384&1&1&0\\
26&82625&0&0&0\\
27&292164&0&0&0\\\bottomrule
\end{tabular}
\end{table}

\begin{proof}
Let $A$ be a facet of $P$ with $24\le m\le27$ faces, and let $G\in\mathcal G_m$ be isomorphic to $\Gamma(A)$. By Lemma~\ref{lem:regular-sets} and~\eqref{eq:R2}, the set $U(A)$ corresponds to a candidate set of size at least $R_2(G)$ in $G$, and $G$ is not $\Gamma(L_{11})$ or $\Gamma(L_{12})$, by Proposition~\ref{prop:small}. An exhaustive search shows that eleven graphs carry such a candidate set, each exactly one (Table~\ref{tab:stage2}). The ten graphs with $24$ vertices, listed in Appendix~\ref{subsec:terminal-data}, are excluded as in Lemma~\ref{lem:stage1}: in each of them the candidate set contains a face which is not a regular pentagon. Together with the $78$ graphs of Lemma~\ref{lem:stage1}, this makes $88$ polyhedra for which enclosures of the normals are used in this section.

The remaining graph $G=(25,23383)$ needs a different argument. It has eighteen vertices of degree five, three of degree six, three of degree seven and one of degree nine, and its only candidate set consists of all eighteen vertices of degree five. So if $G\cong\Gamma(A)$, then all eighteen pentagons of $A$ are regular (Figure~\ref{fig:nonagon}(a)). The nonagon $O$ of $A$ is adjacent to three hexagons and six pentagons, in the cyclic pattern hexagon, pentagon, pentagon, repeated three times. Each of the three hexagons $H$ has, in cyclic order starting from $O$, the neighbours
\[
 O,\ \text{pentagon},\ \text{pentagon},\ \text{heptagon},\ \text{pentagon},\ \text{pentagon}
\]
(Appendix~\ref{subsec:terminal-data}). Since a side shared with a regular pentagon has length $\ell_0$, the side lengths of $H$ in this order are $(x,\ell_0,\ell_0,y,\ell_0,\ell_0)$. Put $a=\cosh x$ and $b=\cosh y$. Formula~\eqref{eq:hexagon}, applied to the sides of lengths $x$ and $y$, gives
\[
 \varphi a=\varphi^2+b,\qquad \varphi b=\varphi^2+a,
\]
so $a=b=\varphi^2/(\varphi-1)=\varphi^3=\cosh2\ell_0$. Thus the sides of $O$ have lengths
\begin{equation}\label{eq:nonagon}
 (2\ell_0,\ell_0,\ell_0,\,2\ell_0,\ell_0,\ell_0,\,2\ell_0,\ell_0,\ell_0)
\end{equation}
in cyclic order. Now, a right-angled polygonal path is determined up to isometry by its sequence of side lengths and the direction of turning (Figure~\ref{fig:nonagon}(b)). The union of a regular right-angled pentagon and its reflection in one of its sides is a right-angled hexagon with side lengths $(2\ell_0,\ell_0,\ell_0,2\ell_0,\ell_0,\ell_0)$. Hence the first six sides of $O$ form a closed hexagon, which is impossible for the boundary of a convex nonagon. So $G$ is not isomorphic to $\Gamma(A)$ either.
\end{proof}

\begin{proof}[Proof of Proposition~\ref{prop:classification}]
Proposition~\ref{prop:upper} leaves $D$, $E$ and facets with $20$ to $27$ faces. By Lemma~\ref{lem:stage1}, those with $20$ to $23$ faces are $\Tt$ and $\Th$, and by Lemma~\ref{lem:stage2} there are none with $24$ to $27$ faces.
\end{proof}
\section{Dissection along internal pentagons}\label{sec:dissection}

In this section $P\subset\HH^4$ is a compact right-angled polytope all of whose facets are isometric to members of $\cA=\{D,E,\Tt,\Th\}$. No bound on the number of facets is assumed. We prove the following theorem and then deduce Theorem~\ref{thm:main}.

\begin{theorem}\label{thm:dissection}
Let $P\subset\HH^4$ be a compact right-angled polytope all of whose facets are isometric to members of $\cA$. Then $P$ can be dissected into $q\ge1$ isometric copies of $Z$ by cutting along $q-1$ hyperplanes, each of which meets the piece being cut in a dodecahedron orthogonal to every facet it meets. Moreover, $f_3(P)=106q+14$.
\end{theorem}

Polytopes with $q>1$ do exist. For example, in the case $q=2$, two copies of the right-angled $120$-cell $Z$ glued along a common facet (Section~\ref{subsec:unions}) form a compact right-angled polytope with $2\cdot119-12=226$ facets, twelve of which are copies of $E$, and cutting along the hyperplane of the common facet recovers the two copies of $Z$.

\subsection{Seams}\label{subsec:seams}
\begin{figure}[ht]
\centering
\begin{tikzpicture}[scale=1.2,line width=0.5pt,line join=round,
  seam/.style={densely dashed,line width=0.8pt}]
\input{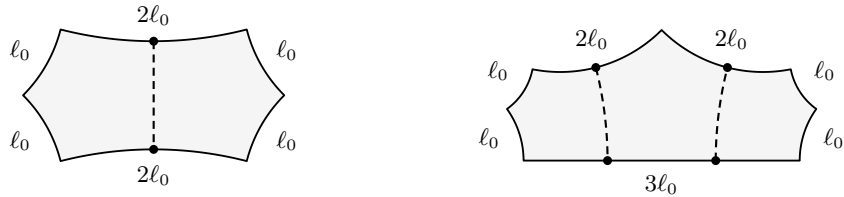}
\begin{scope}[font=\footnotesize]
  \node[above] at (0,0.66) {$2\ell_0$};
  \node[below] at (0,-0.66) {$2\ell_0$};
  \node at (-1.47,0.50) {$\ell_0$};
  \node at (-1.47,-0.50) {$\ell_0$};
  \node at (1.47,0.50) {$\ell_0$};
  \node at (1.47,-0.50) {$\ell_0$};
  \node[below] at (5.6,-0.76) {$3\ell_0$};
  \node at (4.83,0.62) {$2\ell_0$};
  \node at (6.37,0.62) {$2\ell_0$};
  \node at (3.80,0.27) {$\ell_0$};
  \node at (3.70,-0.47) {$\ell_0$};
  \node at (7.40,0.27) {$\ell_0$};
  \node at (7.50,-0.47) {$\ell_0$};
\end{scope}
\end{tikzpicture}
\caption{A hexagonal and a heptagonal face of a member of $\cA$, drawn in the Poincar\'e disk model. The dashed segments are the seams. Their endpoints are the midpoints of the sides of length $2\ell_0$ and the points of trisection of the side of length $3\ell_0$.}\label{fig:seams}
\end{figure}

Let $B$ be a polyhedron in the set $\cA$. The intersection of an internal pentagon of $B$ with a face of $B$ is either empty or a side of the internal pentagon, and in the latter case we call it a \emph{seam} of that face. By the description at the end of Section~\ref{subsec:unions}, the faces of members of $\cA$ are of three kinds (see Figure~\ref{fig:seams}).
\begin{itemize}[leftmargin=*]
\item Pentagons: regular, without seams.
\item Hexagons: unions of two regular pentagons, with side lengths $(2\ell_0,\ell_0,\ell_0,2\ell_0,\ell_0,\ell_0)$ in cyclic order. There is one seam, which joins the midpoints of the two sides of length $2\ell_0$.
\item Heptagons: unions of three regular pentagons, with side lengths $(2\ell_0,\ell_0,\ell_0,3\ell_0,\ell_0,\ell_0,2\ell_0)$ in cyclic order. There are two seams. Each starts at the midpoint of a side of length $2\ell_0$, perpendicular to that side, and ends at one of the two points dividing the side of length $3\ell_0$ into three equal parts.
\end{itemize}
In every case each seam lies in a single face of $B$, its relative interior lies in the relative interior of that face, and its endpoints lie in the relative interiors of edges of $B$.

\begin{lemma}\label{lem:seams}
Every isometry between faces of members of $\cA$ maps seams to seams.
\end{lemma}

\begin{proof}
An isometry preserves side lengths. In a hexagon, the two sides of length $2\ell_0$ and their midpoints are therefore preserved, hence so is the geodesic segment joining them. In a heptagon the sides of length $2\ell_0$ and their midpoints are preserved, and each seam is the segment of the perpendicular at such a midpoint up to the opposite side. So the pair of seams is preserved, although the two seams may be interchanged.
\end{proof}

\subsection{Extending an internal pentagon}
Let $A$ be a facet of $P$ which is not a dodecahedron, let $Q_A$ be an internal pentagon of $A$, and let $\Pi_A\subset H_A$ be the plane in the hyperplane $H_A$ containing $Q_A$. There is a unique hyperplane $H\subset\HH^4$ containing $\Pi_A$ and orthogonal to $H_A$, and it satisfies $H\cap H_A=\Pi_A$ and $H\cap A=Q_A$. The hyperplane $H$ meets the interior of $P$: a point of $\relint Q_A$ lies in $\relint A$, where all other facet inequalities of $P$ are strict, and a short segment in $H$ starting at this point and orthogonal to $H_A$ enters $\inter P$. Consequently, the \emph{section}
\[
 S=P\cap H
\]
of $P$ by $H$ is a compact convex $3$-polyhedron in $H$.

\begin{figure}[ht]
\centering
\begin{tikzpicture}[line width=0.5pt,line join=round,lab/.style={font=\footnotesize}]
\begin{scope}
  \fill[black!5] (-1.560,-2.184) -- (1.560,-2.184) -- (1.560,2.184) -- (-1.560,2.184) -- cycle;
  \fill[black!25] (-0.624,0.000) -- (0.624,0.000) -- (1.010,-1.187) -- (0.000,-1.920) -- (-1.010,-1.187) -- cycle;
  \draw (-0.624,0.000) -- (0.624,0.000) -- (1.010,-1.187) -- (0.000,-1.920) -- (-1.010,-1.187) -- cycle;
  \fill[black!10] (-1.092,-0.421) -- (-0.156,0.421) -- (0.741,0.667) -- (1.092,0.421) -- (0.156,-0.421) -- (-0.741,-0.667) -- cycle;
  \draw[line width=0.8pt] (-1.092,-0.421) -- (-0.156,0.421) -- (0.741,0.667) -- (1.092,0.421) -- (0.156,-0.421) -- (-0.741,-0.667) -- cycle;
  \draw[line width=0.9pt,densely dashed] (-0.624,0.000) -- (0.624,0.000);
  \fill[black!25,fill opacity=0.85] (-0.624,0.000) -- (0.624,0.000) -- (1.010,1.187) -- (0.000,1.920) -- (-1.010,1.187) -- cycle;
  \draw (-0.624,0.000) -- (0.624,0.000) -- (1.010,1.187) -- (0.000,1.920) -- (-1.010,1.187) -- cycle;
  \fill (-0.624,0.000) circle (1.3pt);
  \fill (0.624,0.000) circle (1.3pt);
  \node[lab] at (0.000,-1.131) {$Q_B$};
  \node[lab] at (0.000,1.131) {$Q_C$};
  \node[lab] at (-0.371,-0.369) {$F$};
  \node[lab,anchor=south,inner sep=1pt] at (0.000,0.023) {$e$};
  \node[lab] at (1.365,1.989) {$H$};
  \node[lab,anchor=east] at (-1.365,-0.780) {$B$};
  \node[lab,anchor=east] at (-1.365,0.780) {$C$};
  \node[lab,anchor=south west,inner sep=1pt] at (1.092,0.421) {edge of $P$};
  \node[font=\small] at (0.000,-2.535) {(a)};
\end{scope}
\begin{scope}[xshift=6.6cm,yshift=-0.2cm]
  \fill[black!22] (-0.361,0.498) -- (-0.585,-0.190) -- (-0.000,-0.615) -- (0.585,-0.190) -- (0.361,0.498) -- cycle;
  \fill[black!22] (-0.843,1.161) -- (-0.361,0.498) -- (-0.585,-0.190) -- (-1.365,-0.443) -- (-1.365,0.443) -- cycle;
  \fill[black!10] (-1.365,-0.443) -- (-0.585,-0.190) -- (-0.000,-0.615) -- (-0.000,-1.435) -- (-0.843,-1.161) -- cycle;
  \draw (0.000,2.050) -- (-1.950,0.633) (0.000,2.050) -- (0.000,1.435) (0.000,1.435) -- (-0.843,1.161) (-0.843,1.161) -- (-1.365,0.443) (-0.843,1.161) -- (-0.361,0.498) (-0.361,0.498) -- (-0.585,-0.190) (-1.950,0.633) -- (-1.205,-1.658) (-1.950,0.633) -- (-1.365,0.443) (-1.365,0.443) -- (-1.365,-0.443) (-1.365,-0.443) -- (-0.843,-1.161) (-1.365,-0.443) -- (-0.585,-0.190) (-0.585,-0.190) -- (-0.000,-0.615) (-1.205,-1.658) -- (1.205,-1.658) (-1.205,-1.658) -- (-0.843,-1.161) (-0.843,-1.161) -- (-0.000,-1.435) (-0.000,-1.435) -- (0.843,-1.161) (-0.000,-1.435) -- (-0.000,-0.615) (-0.000,-0.615) -- (0.585,-0.190) (1.205,-1.658) -- (1.950,0.633) (1.205,-1.658) -- (0.843,-1.161) (0.843,-1.161) -- (1.365,-0.443) (1.365,-0.443) -- (1.365,0.443) (1.365,-0.443) -- (0.585,-0.190) (0.585,-0.190) -- (0.361,0.498) (1.950,0.633) -- (0.000,2.050) (1.950,0.633) -- (1.365,0.443) (1.365,0.443) -- (0.843,1.161) (0.843,1.161) -- (0.000,1.435) (0.843,1.161) -- (0.361,0.498) (0.361,0.498) -- (-0.361,0.498);
  \draw[line width=1.3pt] (-0.361,0.498) -- (-0.585,-0.190);
  \fill (-0.585,-0.190) circle (1.5pt);
  \node[lab] at (-0.000,-0.000) {$Q_B$};
  \node[lab] at (-0.912,0.296) {$Q_C$};
  \node[lab] at (-0.563,-0.775) {$Q_{B'}$};
  \node[lab,anchor=east,inner sep=1.5pt] at (-0.473,0.154) {$e$};
  \node[lab,anchor=north,inner sep=2pt] at (-0.585,-0.190) {$p$};
  \node[font=\small] at (0,-2.500) {(b)};
\end{scope}
\end{tikzpicture}
\caption{(a)~The extension step of Lemma~\ref{lem:extension}, drawn schematically in three dimensions: the facets $B$ and $C$, which in $P$ meet at a right angle along the $2$-face $F=B\cap C$ (here a hexagon), are shown on the two sides of $F$. The internal pentagons $Q_B\subset B$ and $Q_C\subset C$ are orthogonal to $F$ and meet it in its seam $e$, so both lie in the hyperplane $H$. In $P$ they meet at a right angle along $e$, as adjacent faces of the section $S$. The endpoints of $e$ are interior points of two edges of $P$. (b)~The section $S=P\cap H$ of Proposition~\ref{prop:section}, a dodecahedron whose faces are the internal pentagons $Q_B$, $B\in\mathcal B$. Adjacent faces $Q_B,Q_C$ meet along the seam $e$ of the $2$-face $B\cap C$, and each vertex $p$ lies in the interior of an edge of $P$, which is contained in exactly three facets of $\mathcal B$.}\label{fig:section}
\end{figure}
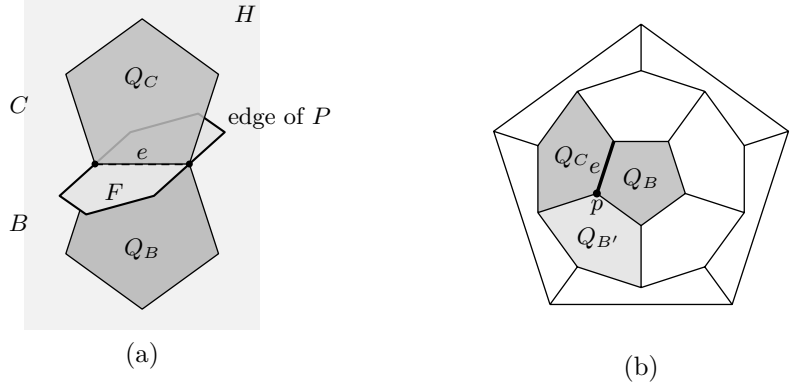

\begin{lemma}\label{lem:extension}
Let $B$ be a facet of $P$ such that $H\perp H_B$ and $H\cap B$ is an internal pentagon $Q_B$ of $B$ (Figure~\ref{fig:section}(a)). Let $e$ be an edge of $Q_B$ and $C$ the neighbour of $B$ across the face of $B$ containing $e$. Then $H\perp H_C$, and $H\cap C$ is an internal pentagon $Q_C$ of $C$, with $Q_B\cap Q_C=e$.
\end{lemma}

\begin{proof}
The edge $e$ is a seam of the $2$-face $F=B\cap C$. By Lemma~\ref{lem:seams}, the edge $e$ is also a seam of $F$ regarded as a face of $C$, so $C$ is not a dodecahedron and $e$ is a side of a unique internal pentagon $Q_C$ of $C$.

Let $n_B,n_C,u$ be unit normals of $H_B,H_C,H$. The right angles give $\langle n_B,n_C\rangle=0$ and, by hypothesis, $\langle u,n_B\rangle=0$. Inside $H_B$ the plane $H\cap H_B$, which contains $Q_B$, is orthogonal to the plane $H_B\cap H_C$ of $F$, because an internal pentagon is orthogonal to every face it meets. As $u$ and $n_C$ are both orthogonal to $n_B$, they are normals within $H_B$ of these two planes, and so $\langle u,n_C\rangle=0$, that is, $H\perp H_C$. Now $H\cap H_C$ and the plane $\Pi_C$ of $Q_C$ are both planes in $H_C$ which contain the line of $e$ and are orthogonal to the plane of $F$. Such a plane is unique, so $H\cap H_C=\Pi_C$, and therefore $H\cap C=\Pi_C\cap C=Q_C$.
\end{proof}

Starting from $A$ and applying Lemma~\ref{lem:extension} repeatedly, the hyperplane $H$ stays fixed throughout, and we obtain a set $\mathcal B$ of facets $B$ of $P$, each with an internal pentagon $Q_B=H\cap B$. A facet with two internal pentagons cannot contribute both, because $H\cap H_B$ is a single plane and the two internal pentagons of a triple union lie in disjoint planes. Thus $Q_B$ is well defined, and $\mathcal B$ is closed under the operation of Lemma~\ref{lem:extension}.

\subsection{The section is a dodecahedron}

\begin{proposition}\label{prop:section}
Let $A$, $Q_A$, $H$, $S$ and $\mathcal B$ be as above. The section $S=P\cap H$ is a dodecahedron. The hyperplane $H$ meets exactly twelve facets of $P$, namely the facets in $\mathcal B$, it is orthogonal to each of them, and it contains no vertex of $P$.
\end{proposition}

\begin{proof}
For $B\in\mathcal B$ we have $S\cap H_B=H\cap B=Q_B$, so $Q_B$ is a $2$-dimensional face of $S$. Let $e$ be an edge of $Q_B$. Its relative interior lies in the relative interior of a $2$-face $B\cap C$ of $P$, so near such a point $S$ is bounded by $H_B$ and $H_C$ only, and the face of $S$ adjacent to $Q_B$ along $e$ is $H\cap C$, which is $Q_C$ by Lemma~\ref{lem:extension}. Thus the faces $Q_B$, $B\in\mathcal B$, form a set of faces of $S$ closed under adjacency (Figure~\ref{fig:section}(b)). Since the face dual graph of a convex polyhedron is connected, these are all the faces of $S$.

All faces of $S$ are therefore regular right-angled pentagons. Two adjacent faces $Q_B,Q_C$ meet at a right angle inside $H$, because $n_B$ and $n_C$ are orthogonal and both tangent to $H$. So $S$ is a compact right-angled polyhedron with only pentagonal faces, and by~\eqref{eq:euler} and the remark after Lemma~\ref{lem:5conn}, $S$ is isometric to $D$. In particular $\mathcal B$ has exactly twelve elements.

It remains to see that no other facet meets $H$. Since $H$ meets $\inter P$, we have $\relint S\subset\inter P$, so a facet of $P$ can only meet $H$ in $\partial S$. A point in the relative interior of a face $Q_B$ of $S$ lies in $\relint B$, hence in no other facet. A point in the relative interior of an edge $Q_B\cap Q_C$ of $S$ lies in the relative interior of the $2$-face $B\cap C$ of $P$, hence only in $B$ and $C$. Finally, a vertex of $S$ is a vertex of an internal pentagon $Q_B$, and hence, by the description in Section~\ref{subsec:seams}, lies in the relative interior of an edge of $P$, which is contained in exactly three facets of $P$, namely the three facets in $\mathcal B$ whose pentagons meet at this vertex. In particular no vertex of $P$ lies in $H$.
\end{proof}

\subsection{Proof of Theorem~\ref{thm:dissection}}
\begin{proof}
For a polytope $P$ as in the theorem let
\[
 c(P)=\#\{\text{facets isometric to }E\}+2\,\#\{\text{facets isometric to }\Tt\text{ or }\Th\}.
\]
We argue by induction on $c(P)$.

Suppose first that $c(P)=0$, so that every facet of $P$ is a dodecahedron. Choose a facet $A$ and move $Z$ by an isometry so that $A$ is a facet of $Z$ and $P$ and $Z$ lie on the same side of $H_A$. If $B$ is the neighbour of $A$ across a pentagon $F$ of $A$, then $H_B$ is the hyperplane containing $F$ and orthogonal to $H_A$, and $B$ is the dodecahedron in $H_B$ which has $F$ as a face and lies on the same side of $H_A$ as $P$. The same description applies to the facet of $Z$ across $F$, so the two coincide. Since $P$ and $Z$ both lie on the side of $H_B$ containing $A$, they lie on the same side of $H_B$. Since $\Gamma(Z)$ is connected, all facets of $Z$ are facets of $P$, with the same inward sides. Hence $P\subseteq Z$. Since $P$ contains all vertices of $Z$, we get $P=Z$.

Now suppose $c(P)>0$. Choose an internal pentagon of some facet and let $H$ be the corresponding hyperplane. By Proposition~\ref{prop:section}, the hyperplane $H$ cuts $P$ into two compact convex polytopes $P_+$ and $P_-$, each having $S=P\cap H$ as a facet. The other facets of $P_\pm$ are the facets of $P$ on the corresponding side of $H$ and the halves of the twelve facets in $\mathcal B$. Every dihedral angle of $P_\pm$ is an old dihedral angle of $P$ or an angle between $S$ and a facet in $\mathcal B$, so $P_\pm$ are right-angled. Cutting $E$ along its internal pentagon gives two copies of $D$, and cutting $\Tt$ or $\Th$ along one of its internal pentagons gives a copy of $D$ and a copy of $E$. So the facets of $P_\pm$ again belong to $\cA$, and since each of the twelve cut facets loses one unit of $c$,
\[
 c(P_+)+c(P_-)=c(P)-12 .
\]
In particular $c(P_\pm)<c(P)$, and the induction hypothesis applies to $P_+$ and $P_-$. This proves the first assertion of Theorem~\ref{thm:dissection}.

For the facet count, note that each cut replaces one polytope by two, which together have the facets of the old one, plus the twelve cut facets counted a second time, plus two copies of $S$:
\[
 f_3(P_+)+f_3(P_-)=f_3(P)+14 .
\]
A dissection into $q$ pieces involves $q-1$ cuts, and every final piece is a copy of $Z$. Summing over all cuts gives $120q=f_3(P)+14(q-1)$, that is, $f_3(P)=106q+14$.
\end{proof}

\subsection{Proof of Theorem~\ref{thm:main}}
\begin{proof}
Let $P\subset\HH^4$ be a compact right-angled polytope with $N\le120$ facets. By Proposition~\ref{prop:classification}, all facets of $P$ belong to $\cA$, and by Theorem~\ref{thm:dissection}, $N=106q+14$ for some $q\ge1$. Hence $q=1$, so $P$ is isometric to $Z$, and $N=120$. Thus no compact right-angled polytope has fewer than $120$ facets, and $Z$ is the only one with exactly $120$.
\end{proof}
\section{Face numbers and volume}\label{sec:consequences}

For a compact right-angled polytope $P\subset\HH^4$ define the \emph{total excess} of its $2$-faces,
\[
 \Delta(P)=\sum_{F}\bigl(k(F)-5\bigr)\ge0,
\]
the sum running over all $2$-faces $F$ of $P$. Thus $\Delta(P)=0$ if and only if all $2$-faces are pentagons, and $\Delta(P)\ge1$ as soon as some facet of $P$ is not a dodecahedron, since such a facet has a face with at least six sides by~\eqref{eq:euler}.

\begin{lemma}\label{lem:fvector}
Let $N=f_3(P)$. Then
\[
 (f_0,f_1,f_2,f_3)(P)=\bigl(5N+\Delta(P),\ 10N+2\Delta(P),\ 6N+\Delta(P),\ N\bigr)
\]
and
\[
 \Vol(P)=\frac{\pi^2}{12}\bigl(f_0(P)-4N+16\bigr)=\frac{\pi^2}{12}\bigl(N+\Delta(P)+16\bigr).
\]
\end{lemma}

\begin{proof}
Since $P$ is simple, every vertex lies on four edges and on six $2$-faces. Hence $f_1=2f_0$ and $\sum_F k(F)=6f_0$, and Euler's relation $f_0-f_1+f_2-f_3=0$ gives $f_2=f_0+N$. Therefore $\Delta(P)=6f_0-5f_2=f_0-5N$, which gives the face vector.

For the volume we use the Gauss--Bonnet formula for the hyperbolic orbifold $\HH^4/W$, where $W$ is the reflection group of $P$: $\Vol(P)=\tfrac{4\pi^2}{3}\chi(W)$, where $\chi(W)=\sum_F(-1)^{\dim F}/|W_F|$ is the orbifold Euler characteristic, the sum running over all faces $F$ of $P$ including $P$ itself, and $W_F$ is the stabilizer of a point in the relative interior of $F$. Since $P$ is right-angled, $|W_F|=2^{4-\dim F}$, and
\[
 \chi(W)=1-\frac{f_3}{2}+\frac{f_2}{4}-\frac{f_1}{8}+\frac{f_0}{16}=\frac{f_0-4N+16}{16}
\]
by the relations above. This is the formula of Kellerhals and Perren~\cite[Proposition~3.2]{KP11}.
\end{proof}

\begin{proof}[Proof of Corollary~\ref{cor:volume}]
By Lemma~\ref{lem:fvector}, each of $f_0,f_1,f_2$ and $\Vol$ is an increasing function of $N$ and of $\Delta(P)\ge0$. By Theorem~\ref{thm:main}, $N\ge120$. This gives the four inequalities, with equality only if $N=120$, that is, only if $P\cong Z$. Conversely, the $120$-cell $Z$ has $N=120$ and $\Delta(Z)=0$, and so attains all four bounds. In particular, $\Vol(Z)=136\pi^2/12=34\pi^2/3$.
\end{proof}

\begin{proof}[Proof of Corollary~\ref{cor:kolpakov}]
If $P$ has a facet that is not isometric to $D$, then $P\not\cong Z$, so $N\ge121$ by Theorem~\ref{thm:main}, and $\Delta(P)\ge1$. Lemma~\ref{lem:fvector} gives $f_0(P)\ge5\cdot121+1=606$ and $\Vol(P)\ge\pi^2(121+1+16)/12=23\pi^2/2$.
\end{proof}

\subsection*{Questions}
Two copies of $Z$ glued along a facet give a compact right-angled polytope with $226$ facets (Section~\ref{sec:dissection}). By Theorems~\ref{thm:main} and~\ref{thm:dissection}, the second smallest number of facets of a compact right-angled polytope in $\HH^4$ lies between $121$ and $226$, and it is attained, if it is below $226$, only by polytopes having a facet that is not isometric to any member of $\cA$.

\begin{question}
What is the second smallest number of facets of a compact right-angled polytope in $\HH^4$?
\end{question}

In dimension three, the L\"obell polyhedron $L_6$ has the second smallest volume among compact right-angled polyhedra, by a theorem of Inoue~\cite{Ino08}.

\begin{question}
Which compact right-angled polytope in $\HH^4$ has the second smallest volume?
\end{question}
\appendix

\section{The computations}\label{app:computation}

This appendix describes the computations used in Sections~\ref{sec:small}--\ref{sec:classification} and explains why their output proves what is claimed. All exclusions are based either on rigorous enclosures, computed and compared in exact rational arithmetic, or on exact algebraic computation. Floating-point numbers are only used to find approximate data, which are then verified.

\subsection{Catalogues}\label{subsec:catalogues}
For $12\le m\le27$ let $\mathcal G_m$ be the list of $5$-connected triangulations of the $2$-sphere with $m$ vertices in the supplementary material, each given as a rotation system (the cyclically ordered list of neighbours of every vertex). The lists were produced with \texttt{plantri}~\cite{BM07}. We write $(m,i)$ for the $i$-th entry of $\mathcal G_m$, counting from $0$.

Each rotation system is checked to define a triangulation of the sphere which is flag and has no squares, and the members of each list are checked to be pairwise nonisomorphic, by comparing canonical forms. The cardinalities are
\begin{center}
\begin{tabular}{lrrrrrrrr}\toprule
$m$ & 12&13&14&15&16&17&18&19\\
$|\mathcal G_m|$ &1&0&1&1&3&4&12&23\\\midrule
$m$ & 20&21&22&23&24&25&26&27\\
$|\mathcal G_m|$ &71&187&627&1970&6833&23384&82625&292164\\\bottomrule
\end{tabular}
\end{center}
These agree with the numbers of $5$-connected triangulations computed by Brinkmann and McKay~\cite[Table~4, column $t_{5\text{-}5}$, in the expanded version]{BM07}, which count isomorphism classes with orientation-reversing isomorphisms allowed, that is, isomorphism classes of graphs. Hence each $\mathcal G_m$ contains exactly one representative of every isomorphism class, and by Lemma~\ref{lem:5conn} it contains the graph of every compact right-angled polyhedron with $m$ faces. In total there are $407906$ graphs.

\subsection{Verified realizations}\label{subsec:realizations}
Let $G\in\mathcal G_m$ and let $B$ be a compact right-angled polyhedron with $\Gamma(B)=G$, the faces of $B$ being labelled by the vertices of $G$. Its outward unit normals $n_0,\dots,n_{m-1}\in\RR^{3,1}$ satisfy
\begin{equation}\label{eq:system}
 \langle n_a,n_a\rangle=1,\qquad \langle n_a,n_b\rangle=0\quad(ab\in E(G)).
\end{equation}
Choose a triangle $abc$ of $G$. After an isometry we may assume $n_a,n_b,n_c$ are the first three standard basis vectors. The remaining $m-3$ normals have $M=4m-12$ coordinates $x\in\RR^M$, and~\eqref{eq:system} becomes a system $\Phi(x)=0$ of $(m-3)+(3m-9)=M$ polynomial equations of degree at most two with integer coefficients.

\begin{proposition}\label{prop:contraction}
Let $x_0\in\mathbb Q^M$, let $C$ be a rational $M\times M$ matrix and $\varepsilon>0$ rational. For $1\le i\le M$ put
\[
 \eta_i=\bigl|(C\Phi(x_0))_i\bigr|,\qquad \kappa_i=\sum_j\bigl|(I-C\,D\Phi(x_0))_{ij}\bigr|,\qquad \lambda_i=\sum_j|C_{ij}|\,L_j ,
\]
where $L_j=0$ if the $j$-th equation is linear and $L_j=8$ otherwise, and suppose that $\eta_i+\varepsilon\kappa_i+\varepsilon^2\lambda_i<\varepsilon$ for all $i$. Then $\Phi$ has a unique zero $x_*$ with $\|x_*-x_0\|_\infty\le\varepsilon$.
\end{proposition}

\begin{proof}
Every component $\Phi_j$ is either linear, or of the form $\langle y,y\rangle-1$ or $\langle y,y'\rangle$ for $4$-vectors $y,y'$ of unknowns. Hence $\|D\Phi_j(x)-D\Phi_j(x_0)\|_1\le L_j\|x-x_0\|_\infty$. For $\Psi(x)=x-C\Phi(x)$ and $\|x-x_0\|_\infty\le\varepsilon$, Taylor's formula gives
\[
 |(\Psi(x)-x_0)_i|\le\eta_i+\kappa_i\varepsilon+\tfrac12\lambda_i\varepsilon^2<\varepsilon,
\]
and the $i$-th row of $D\Psi(x)$ has $\ell^1$-norm at most $\kappa_i+\lambda_i\varepsilon<1$. So $\Psi$ is a contraction of the cube $\|x-x_0\|_\infty\le\varepsilon$ into itself and has a unique fixed point $x_*$ there. Since $\|I-C\,D\Phi(x_0)\|_\infty<1$, the matrix $C$ is invertible, and $\Psi(x_*)=x_*$ means $\Phi(x_*)=0$.
\end{proof}

This is a standard device in validated numerics (see~\cite{Tuc11}). In practice $x_0$ is obtained by Newton's method in floating point, $C$ is a rational approximation of $D\Phi(x_0)^{-1}$, and $x_0,C,\varepsilon$ are stored in the supplementary material. The inequalities of Proposition~\ref{prop:contraction} are then verified in exact rational arithmetic. Such data were produced for $195$ polyhedra.

\subsection{Compactness and combinatorial type}
The zero $x_*$ of Proposition~\ref{prop:contraction} need not a priori come from a polyhedron with the right combinatorics. This is checked as follows. Using interval arithmetic over the whole cube $\|x-x_0\|_\infty\le\varepsilon$, we verify:
\begin{enumerate}
\item $\langle n_a,n_b\rangle<-1$ for every pair of nonadjacent vertices $a,b$ of $G$;
\item for every triangle $abc$ of $G$, the last coordinate of the Lorentzian cross product $v_{abc}$ of $n_a,n_b,n_c$ (a vector orthogonal to all three, given by $3\times3$ minors) does not vanish, and after changing the sign of $v_{abc}$, if necessary, to make this coordinate positive, $\langle v_{abc},n_j\rangle<0$ for all $j\notin\{a,b,c\}$.
\end{enumerate}

\begin{lemma}\label{lem:type}
If~(i) and~(ii) hold, then the halfspaces $\{\langle x,n_j\rangle\le0\}$ cut out a compact right-angled polyhedron $B_*\subset\HH^3$ whose face dual graph is $G$, with the given labels.
\end{lemma}

\begin{proof}
At $x_*$ the vectors $n_a,n_b,n_c$ of a triangle are orthonormal, so $v_{abc}$, suitably normalized, is a point of $\HH^3$ lying on the three planes $a,b,c$ and strictly inside all other halfspaces. Work in the Klein model, and let $B_\cap$ be the intersection of the $m$ halfspaces and $B_{\rm hull}$ the convex hull of the points $v_{abc}$. Near $v_{abc}$, $B_\cap$ coincides with the simplicial cone cut out by $a,b,c$, and its three edges contain the points $v_{abd}$, $v_{acd'}$, $v_{bcd''}$ of the three triangles adjacent to $abc$. Hence $B_{\rm hull}$ has the same tangent cone at $v_{abc}$ as $B_\cap$. Every facet of $B_{\rm hull}$ contains some $v_{abc}$ and is therefore contained in one of the planes $a,b,c$. So $B_{\rm hull}$ is an intersection of some of the $m$ halfspaces and contains $B_\cap$. As $B_{\rm hull}\subseteq B_\cap$ by convexity, $B_{\rm hull}=B_\cap=B_*$. Its vertices are the points $v_{abc}$, its faces correspond to the vertices of $G$, two faces are adjacent exactly when the corresponding vertices of $G$ span an edge, and adjacent faces are orthogonal by~\eqref{eq:system}. (Condition~(i) is not needed for this argument. It is a consistency check, which also confirms~(F2) for $B_*$.)
\end{proof}

By Theorem~\ref{thm:andreev}, the polyhedron $B_*$ is therefore the right-angled realization of $G$.

\subsection{Edge lengths}\label{subsec:edges}
Let $ab$ be an edge of $G$, lying in the triangles $abc$ and $abd$. The corresponding edge of $B_*$ is the segment of the line $H_a\cap H_b$ between the planes $H_c$ and $H_d$, both orthogonal to this line, and its length $\ell_{ab}$ is the distance between these two ultraparallel planes. Hence
\begin{equation}\label{eq:edge-length}
 \cosh\ell_{ab}=-\langle n_c,n_d\rangle ,
\end{equation}
and evaluating the right-hand side in interval arithmetic gives a rational enclosure of $\cosh\ell_{ab}$. To show that a face is not a regular pentagon, we exhibit a side whose enclosure is disjoint from a rational enclosure of $\varphi=\cosh\ell_0$. The latter is obtained from rational bounds for $\sqrt5$ verified by squaring. Overlapping enclosures are never interpreted as equality.

For the L\"obell polyhedra exact values are available from~\eqref{eq:lobell}: the edges of $L_6$ have $\cosh^2\ell\in\{2,3,4\}$, and those of $L_7$ have $\cosh^2\ell\in\{\tau^2,\,\tau/(\tau-1),\,\tau+1\}$, where $\tau=1+2\cos(2\pi/7)$ is the root in $(2,3)$ of $\tau^3-2\tau^2-\tau+1$. Several exclusions for these two sources are equality cases of strict inequalities and must be decided exactly. For example, two edges with $\sinh\ell=1$ (the value for type $\ell_{6,2}$) cannot meet at a vertex, by Lemma~\ref{lem:short}. Also, a base edge $e=zz'$ of $L_6$, or an edge $e=zz'$ of type $\ell_{7,3}$ of $L_7$, satisfies $\cosh\ell(e)=\cosh\ell(a)\cosh\ell(a')$ exactly, where $a,a'$ are the edges of type $\ell_{n,2}$ at $z,z'$, so that Lemma~\ref{lem:placement}(iii) applies only because its inequality is strict (compare the proof of Lemma~\ref{lem:lobell}). More generally, for the source $L_6$ most placements are excluded without using normals, by propagating the exact side lengths through the pentagon relations~\eqref{eq:pentagon} inside the target until a contradiction appears.

\subsection{Placements and their compatibility}\label{subsec:matching}
We now describe the finite checks behind Lemmas~\ref{lem:single} and~\ref{lem:simultaneous}. Fix the source polyhedron $A$ and a face $Q$ of $A$ with vertices $z_0,\dots,z_{k-1}$ and sides $e_j=z_{j-1}z_j$. For every target $G\in\mathcal G_{m'}$ with $m'$ in the range given by~\eqref{eq:finite}, every face $R$ of the target with $k(R)=k$ and defect below the cutoff, and each of the $2k$ placements $\iota$ along $R$, the following tests are applied, and a placement is \emph{excluded} if one of them fails.
\begin{enumerate}
\item \emph{Forced nonpentagons.} If the enclosures show that $\cosh\ell(e_j)\ge\cosh\ell(a)\cosh\ell(a')$ for some edges $a$ of $A$ at $z_{j-1}$ and $a'$ of $A$ at $z_j$, then the face of the target across $\iota(e_j)$ must not be a pentagon (Lemma~\ref{lem:placement}(iii)).
\item \emph{Edge orbits.} Edges of the target in one orbit of $\Aut(G)$ have equal length (Theorem~\ref{thm:andreev}). If two sides of $Q$ with disjoint enclosures are mapped into the same orbit, the placement is excluded. This test only needs the source enclosures.
\item \emph{Lengths.} For each $j$, the enclosures of $\cosh\ell(e_j)$ and of $\cosh\ell(\iota(e_j))$ must intersect (Lemma~\ref{lem:placement}(i)).
\item \emph{Transverse edges.} For each $j$, with $t(z_j)$ the edge of the target at $\iota(z_j)$ not on $R$, the inequalities $\sinh\ell(a)\sinh\ell(t(z_j))>1$, for the edges $a$ of $A$ at $z_j$, must not be refuted by the enclosures (Lemma~\ref{lem:placement}(ii)).
\end{enumerate}
Every realized placement passes all four tests, so the set of \emph{surviving} placements contains all realized ones. (A surviving placement need not be realized.) Tests~(i) and~(ii) need no data for the target, and only targets passing them require verified realizations. Placements that differ by an automorphism of the target fixing $R$ have the same lengths and defects, so they may be identified, and when comparing two placements we use a fixed representative in each class.

For the compatibility statements we consider surviving placements of $A$ on itself. Two such placements, on faces $Q$ and $Q'$ of $A$, are joined by an edge of the \emph{obstruction graph} $\Omega_A$ if $Q=Q'$, or if $Q$ and $Q'$ share an edge and, at one of its endpoints $z$, the enclosures of the lengths of the two transverse edges $t(z)$ are disjoint. By Lemma~\ref{lem:placement}(iv), the placements realized by the neighbours of $A$ that are isometric to $A$ form an independent set in $\Omega_A$. Lemma~\ref{lem:simultaneous}(i) and~(ii) are obtained from explicit subgraphs of $\Omega_A$ recorded in Appendix~\ref{subsec:simultaneous-data}.

\subsection{A worked example}\label{subsec:example}
To illustrate Section~\ref{subsec:candidates}, consider the entry $G=(20,16)$ of $\mathcal G_{20}$. It has twelve vertices of degree five, so $R_1(G)=\ceil{(6+240+48-240)/9}=6$, and the search finds exactly one candidate set,
\[
 U=\{4,5,8,9,12,13,16,17\},
\]
in the labelling of the supplementary material. The face $4$ and its neighbour $3$ share an edge, and the edge $\{3,4\}$ of $G$ lies in the triangles $\{3,4,11\}$ and $\{3,4,18\}$. The verified realization of $G$ and~\eqref{eq:edge-length} give, after outward rounding,
\[
 \cosh\ell_{34}=-\langle n_{11},n_{18}\rangle\in\left[\frac{1274309}{10^6},\frac{1274310}{10^6}\right],
 \qquad
 \varphi\in\left[\frac{1618033}{10^6},\frac{1618034}{10^6}\right].
\]
The intervals are disjoint, so face $4$ is not a regular pentagon and $U$ is not the set of regular pentagons of the right-angled realization of $G$. Since $U$ is the only candidate, no facet of a polytope with at most $120$ facets has $G$ as its graph. (Of course, $G$ is the graph of a compact right-angled polyhedron. It simply cannot occur as such a facet.)

\subsection{Supplementary material and independent checks}\label{subsec:supplement}
The supplementary material, available at \url{https://github.com/Lizi-Guo/120cell-minimality}, contains the catalogues $\mathcal G_m$, the rational data $(x_0,C,\varepsilon)$ for the $195$ verified realizations, the candidate sets and placements, and programs that verify every finite statement used in the paper. Its README describes the input formats and how to rerun the computations, and a complete run takes a few minutes on a laptop.

The candidate search of Section~\ref{sec:classification} was implemented twice, once as a recursive subset enumeration and once by propagation of Boolean constraints, and the two programs produce the same labelled candidate sets. In addition, all combinatorial counts in Sections~\ref{sec:small}--\ref{sec:classification} (the catalogue sizes, Tables~\ref{tab:stage1} and~\ref{tab:stage2}, the table in the proof of Lemma~\ref{lem:clusters}, the defects of $E$, $\Tt$, $\Th$, and the minima used in Lemma~\ref{lem:simultaneous}(iv)) were recomputed by a separately written program, starting from a catalogue regenerated with \texttt{plantri}, and the graphs it found were matched with the stored catalogue entries through canonical forms. That program also recomputed in floating point, by Newton's method, all realizations used in the proof, and reproduced the conclusions of Lemmas~\ref{lem:small-facts}, \ref{lem:single}, \ref{lem:simultaneous}, \ref{lem:stage1} and~\ref{lem:stage2}: the same graphs are excluded and the same placements survive. Apart from the exact comparisons for $L_6$ and $L_7$ described in Section~\ref{subsec:edges}, every exclusion rests on a strict inequality that holds with a margin of at least $5\times10^{-3}$ when both sides are written in terms of $\cosh^2$ of edge lengths, and usually with much more. The verified enclosures are far narrower than this.
\section{Data}\label{app:data}

Labels of faces refer to the rotation systems of the catalogue entries in the supplementary material. For the types with at most $17$ faces they are reproduced in Appendix~\ref{subsec:small-data}.

\subsection{The polyhedra with at most seventeen faces}\label{subsec:small-data}
The ten combinatorial types with at most $17$ faces (Table~\ref{tab:small}) are listed below as rotation systems. For each face $Q$ of the polyhedron $B$ we give its number of sides $k=k(Q)$, its neighbours in cyclic order, and the defect $\dd=\dd(B,Q)=m(B)-12+\dd_{\ge2}(B,Q)$ of~\eqref{eq:defect}. The last table lists the orbits of the faces under the automorphism group of $\Gamma$, which by Theorem~\ref{thm:andreev} is the isometry group. These are the face orbits of Table~\ref{tab:same-type} and of the last paragraph of Appendix~\ref{subsec:single-data}. The labels are those of the catalogue entries $(m,i)$ of Appendix~\ref{subsec:catalogues}, and they are the labels used in the rest of this appendix. These data give Lemma~\ref{lem:small-facts}(ii) directly, namely $\dd(B,Q)\ge7$ at every face of each of the four types with $17$ faces, as well as the defects $10$, $8$ and $7$ of $E$ used in Section~\ref{sec:identity} and in the proof of Lemma~\ref{lem:firstcost}. They also show that the four hexagons of $F_{16b}$ are pairwise disjoint and that every pentagon of $F_{16b}$ meets exactly two of them (Section~\ref{subsec:small-list}).

\begin{center}\footnotesize\setlength{\tabcolsep}{3.5pt}
\begin{tabular}[t]{rrlr}\toprule
\multicolumn{4}{l}{$D=(12,0)$}\\\midrule
$Q$ & $k$ & neighbours & $\dd$\\\midrule
0 & 5 & 2,6,5,4,3 & 0\\
1 & 5 & 7,8,9,10,11 & 0\\
2 & 5 & 0,3,7,11,6 & 0\\
3 & 5 & 0,4,8,7,2 & 0\\
4 & 5 & 0,5,9,8,3 & 0\\
5 & 5 & 0,6,10,9,4 & 0\\
6 & 5 & 0,2,11,10,5 & 0\\
7 & 5 & 1,11,2,3,8 & 0\\
8 & 5 & 1,7,3,4,9 & 0\\
9 & 5 & 1,8,4,5,10 & 0\\
10 & 5 & 1,9,5,6,11 & 0\\
11 & 5 & 1,10,6,2,7 & 0\\\bottomrule
\end{tabular}
\hfill
\begin{tabular}[t]{rrlr}\toprule
\multicolumn{4}{l}{$L_6=(14,0)$}\\\midrule
$Q$ & $k$ & neighbours & $\dd$\\\midrule
0 & 6 & 2,7,6,5,4,3 & 3\\
1 & 6 & 8,9,10,11,12,13 & 3\\
2 & 5 & 0,3,8,13,7 & 3\\
3 & 5 & 0,4,9,8,2 & 3\\
4 & 5 & 0,5,10,9,3 & 3\\
5 & 5 & 0,6,11,10,4 & 3\\
6 & 5 & 0,7,12,11,5 & 3\\
7 & 5 & 0,2,13,12,6 & 3\\
8 & 5 & 1,13,2,3,9 & 3\\
9 & 5 & 1,8,3,4,10 & 3\\
10 & 5 & 1,9,4,5,11 & 3\\
11 & 5 & 1,10,5,6,12 & 3\\
12 & 5 & 1,11,6,7,13 & 3\\
13 & 5 & 1,12,7,2,8 & 3\\\bottomrule
\end{tabular}
\hfill
\begin{tabular}[t]{rrlr}\toprule
\multicolumn{4}{l}{$F_{15}=(15,0)$}\\\midrule
$Q$ & $k$ & neighbours & $\dd$\\\midrule
0 & 5 & 2,7,6,5,14 & 4\\
1 & 6 & 8,9,10,11,12,13 & 5\\
2 & 6 & 0,14,3,8,13,7 & 5\\
3 & 5 & 2,14,4,9,8 & 5\\
4 & 5 & 3,14,5,10,9 & 5\\
5 & 6 & 0,6,11,10,4,14 & 5\\
6 & 5 & 0,7,12,11,5 & 5\\
7 & 5 & 0,2,13,12,6 & 5\\
8 & 5 & 1,13,2,3,9 & 4\\
9 & 5 & 1,8,3,4,10 & 5\\
10 & 5 & 1,9,4,5,11 & 4\\
11 & 5 & 1,10,5,6,12 & 4\\
12 & 5 & 1,11,6,7,13 & 5\\
13 & 5 & 1,12,7,2,8 & 4\\
14 & 5 & 0,5,4,3,2 & 4\\\bottomrule
\end{tabular}
\hfill
\begin{tabular}[t]{rrlr}\toprule
\multicolumn{4}{l}{$L_7=(16,0)$}\\\midrule
$Q$ & $k$ & neighbours & $\dd$\\\midrule
0 & 7 & 2,8,7,6,5,4,3 & 6\\
1 & 7 & 9,10,11,12,13,14,15 & 6\\
2 & 5 & 0,3,9,15,8 & 6\\
3 & 5 & 0,4,10,9,2 & 6\\
4 & 5 & 0,5,11,10,3 & 6\\
5 & 5 & 0,6,12,11,4 & 6\\
6 & 5 & 0,7,13,12,5 & 6\\
7 & 5 & 0,8,14,13,6 & 6\\
8 & 5 & 0,2,15,14,7 & 6\\
9 & 5 & 1,15,2,3,10 & 6\\
10 & 5 & 1,9,3,4,11 & 6\\
11 & 5 & 1,10,4,5,12 & 6\\
12 & 5 & 1,11,5,6,13 & 6\\
13 & 5 & 1,12,6,7,14 & 6\\
14 & 5 & 1,13,7,8,15 & 6\\
15 & 5 & 1,14,8,2,9 & 6\\\bottomrule
\end{tabular}
\end{center}

\begin{center}\footnotesize\setlength{\tabcolsep}{3.5pt}
\begin{tabular}[t]{rrlr}\toprule
\multicolumn{4}{l}{$F_{16a}=(16,1)$}\\\midrule
$Q$ & $k$ & neighbours & $\dd$\\\midrule
0 & 5 & 2,7,6,5,14 & 6\\
1 & 5 & 8,9,10,11,15 & 6\\
2 & 6 & 0,14,3,8,13,7 & 6\\
3 & 5 & 2,14,4,9,8 & 6\\
4 & 5 & 3,14,5,10,9 & 7\\
5 & 6 & 0,6,11,10,4,14 & 6\\
6 & 5 & 0,7,12,11,5 & 6\\
7 & 5 & 0,2,13,12,6 & 7\\
8 & 6 & 1,15,13,2,3,9 & 6\\
9 & 5 & 1,8,3,4,10 & 7\\
10 & 5 & 1,9,4,5,11 & 6\\
11 & 6 & 1,10,5,6,12,15 & 6\\
12 & 5 & 6,7,13,15,11 & 7\\
13 & 5 & 2,8,15,12,7 & 6\\
14 & 5 & 0,5,4,3,2 & 6\\
15 & 5 & 1,11,12,13,8 & 6\\\bottomrule
\end{tabular}
\hfill
\begin{tabular}[t]{rrlr}\toprule
\multicolumn{4}{l}{$F_{16b}=(16,2)$}\\\midrule
$Q$ & $k$ & neighbours & $\dd$\\\midrule
0 & 5 & 2,7,6,5,14 & 6\\
1 & 5 & 9,10,11,12,15 & 6\\
2 & 6 & 0,14,3,8,13,7 & 7\\
3 & 5 & 2,14,4,9,8 & 6\\
4 & 5 & 3,14,5,10,9 & 6\\
5 & 6 & 0,6,11,10,4,14 & 7\\
6 & 5 & 0,7,12,11,5 & 6\\
7 & 5 & 0,2,13,12,6 & 6\\
8 & 5 & 2,3,9,15,13 & 6\\
9 & 6 & 1,15,8,3,4,10 & 7\\
10 & 5 & 1,9,4,5,11 & 6\\
11 & 5 & 1,10,5,6,12 & 6\\
12 & 6 & 1,11,6,7,13,15 & 7\\
13 & 5 & 2,8,15,12,7 & 6\\
14 & 5 & 0,5,4,3,2 & 6\\
15 & 5 & 1,12,13,8,9 & 6\\\bottomrule
\end{tabular}
\hfill
\begin{tabular}[t]{rrlr}\toprule
\multicolumn{4}{l}{$F_{17a}=(17,0)$}\\\midrule
$Q$ & $k$ & neighbours & $\dd$\\\midrule
0 & 5 & 2,8,7,6,16 & 7\\
1 & 7 & 9,10,11,12,13,14,15 & 8\\
2 & 6 & 0,16,3,9,15,8 & 8\\
3 & 5 & 2,16,4,10,9 & 8\\
4 & 5 & 3,16,5,11,10 & 9\\
5 & 5 & 4,16,6,12,11 & 8\\
6 & 6 & 0,7,13,12,5,16 & 8\\
7 & 5 & 0,8,14,13,6 & 9\\
8 & 5 & 0,2,15,14,7 & 9\\
9 & 5 & 1,15,2,3,10 & 7\\
10 & 5 & 1,9,3,4,11 & 8\\
11 & 5 & 1,10,4,5,12 & 8\\
12 & 5 & 1,11,5,6,13 & 7\\
13 & 5 & 1,12,6,7,14 & 7\\
14 & 5 & 1,13,7,8,15 & 8\\
15 & 5 & 1,14,8,2,9 & 7\\
16 & 6 & 0,6,5,4,3,2 & 7\\\bottomrule
\end{tabular}
\hfill
\begin{tabular}[t]{rrlr}\toprule
\multicolumn{4}{l}{$F_{17b}=(17,1)$}\\\midrule
$Q$ & $k$ & neighbours & $\dd$\\\midrule
0 & 5 & 5,14,16,7,6 & 8\\
1 & 5 & 8,9,10,11,15 & 8\\
2 & 5 & 3,8,13,16,14 & 7\\
3 & 5 & 2,14,4,9,8 & 8\\
4 & 5 & 3,14,5,10,9 & 8\\
5 & 6 & 0,6,11,10,4,14 & 7\\
6 & 5 & 0,7,12,11,5 & 8\\
7 & 5 & 0,16,13,12,6 & 9\\
8 & 6 & 1,15,13,2,3,9 & 8\\
9 & 5 & 1,8,3,4,10 & 9\\
10 & 5 & 1,9,4,5,11 & 8\\
11 & 6 & 1,10,5,6,12,15 & 8\\
12 & 5 & 6,7,13,15,11 & 8\\
13 & 6 & 2,8,15,12,7,16 & 8\\
14 & 6 & 0,5,4,3,2,16 & 8\\
15 & 5 & 1,11,12,13,8 & 7\\
16 & 5 & 0,14,2,13,7 & 8\\\bottomrule
\end{tabular}
\end{center}

\begin{center}\footnotesize\setlength{\tabcolsep}{3.5pt}
\begin{tabular}[t]{rrlr}\toprule
\multicolumn{4}{l}{$F_{17c}=(17,2)$}\\\midrule
$Q$ & $k$ & neighbours & $\dd$\\\midrule
0 & 5 & 5,14,16,7,6 & 8\\
1 & 5 & 8,9,10,11,15 & 8\\
2 & 5 & 3,8,13,7,16 & 7\\
3 & 6 & 2,16,14,4,9,8 & 8\\
4 & 5 & 3,14,5,10,9 & 8\\
5 & 6 & 0,6,11,10,4,14 & 8\\
6 & 5 & 0,7,12,11,5 & 7\\
7 & 6 & 0,16,2,13,12,6 & 9\\
8 & 6 & 1,15,13,2,3,9 & 8\\
9 & 5 & 1,8,3,4,10 & 8\\
10 & 5 & 1,9,4,5,11 & 8\\
11 & 6 & 1,10,5,6,12,15 & 8\\
12 & 5 & 6,7,13,15,11 & 8\\
13 & 5 & 2,8,15,12,7 & 8\\
14 & 5 & 0,5,4,3,16 & 8\\
15 & 5 & 1,11,12,13,8 & 8\\
16 & 5 & 0,14,3,2,7 & 8\\\bottomrule
\end{tabular}
\hfill
\begin{tabular}[t]{rrlr}\toprule
\multicolumn{4}{l}{$E=(17,3)$}\\\midrule
$Q$ & $k$ & neighbours & $\dd$\\\midrule
0 & 5 & 6,7,8,9,10 & 10\\
1 & 6 & 2,6,10,5,12,16 & 7\\
2 & 6 & 1,16,15,3,7,6 & 7\\
3 & 6 & 2,15,14,4,8,7 & 7\\
4 & 6 & 3,14,13,5,9,8 & 7\\
5 & 6 & 1,10,9,4,13,12 & 7\\
6 & 5 & 0,10,1,2,7 & 8\\
7 & 5 & 0,6,2,3,8 & 8\\
8 & 5 & 0,7,3,4,9 & 8\\
9 & 5 & 0,8,4,5,10 & 8\\
10 & 5 & 0,9,5,1,6 & 8\\
11 & 5 & 12,13,14,15,16 & 10\\
12 & 5 & 1,5,13,11,16 & 8\\
13 & 5 & 4,14,11,12,5 & 8\\
14 & 5 & 3,15,11,13,4 & 8\\
15 & 5 & 2,16,11,14,3 & 8\\
16 & 5 & 1,12,11,15,2 & 8\\\bottomrule
\end{tabular}
\hfill
\begin{tabular}[t]{l>{\raggedright\arraybackslash}p{5.6cm}}\toprule
\multicolumn{2}{l}{orbits of faces under $\Aut(\Gamma)$}\\\midrule
$D$ & $\{0,1,2,3,4,5,6,7,8,9,10,11\}$\\
$L_6$ & $\{0,1\}$, $\{2,3,4,5,6,7,8,9,10,11,12,13\}$\\
$F_{15}$ & $\{0,8,10,11,13,14\}$, $\{1,2,5\}$, $\{3,4,6,7,9,12\}$\\
$L_7$ & $\{0,1\}$, $\{2,3,4,5,6,7,8,9,10,11,12,13,14,15\}$\\
$F_{16a}$ & $\{0,1,14,15\}$, $\{2,5,8,11\}$, $\{3,6,10,13\}$, $\{4,7,9,12\}$\\
$F_{16b}$ & $\{0,1,3,4,6,7,8,10,11,13,14,15\}$, $\{2,5,9,12\}$\\
$F_{17a}$ & $\{0\}$, $\{1\}$, $\{2,6\}$, $\{3,5\}$, $\{4\}$, $\{7,8\}$, $\{9,12\}$, $\{10,11\}$, $\{13,15\}$, $\{14\}$, $\{16\}$\\
$F_{17b}$ & $\{0,4,6,10\}$, $\{1,3,12,16\}$, $\{2,15\}$, $\{5\}$, $\{7,9\}$, $\{8,13\}$, $\{11,14\}$\\
$F_{17c}$ & $\{0,12,13,16\}$, $\{1,4\}$, $\{2,6\}$, $\{3,5,8,11\}$, $\{7\}$, $\{9,10\}$, $\{14,15\}$\\
$E$ & $\{0,11\}$, $\{1,2,3,4,5\}$, $\{6,7,8,9,10,12,13,14,15,16\}$\\\bottomrule
\end{tabular}
\end{center}

\subsection{Single neighbours}\label{subsec:single-data}
Table~\ref{tab:coverage} summarizes the searches in the proof of Lemma~\ref{lem:single}(i) and~(ii). For $L_6$ the target range contains all $116$ types with at most $20$ faces. The placements that survive the exact propagation of Appendix~\ref{subsec:edges} involve $14$ targets and are then excluded with verified realizations. For $L_7$ two placements, on the targets $(17,1)$ and $(17,2)$, survive the preliminary tests and are excluded in the same way. For $F_{15}$, $F_{16a}$ and $F_{16b}$, verified realizations of $117$ targets are used.

\begin{table}[ht]
\centering
\caption{Searches for single neighbours. The third column gives the number of placements, with one target face in each orbit of target symmetries, that pass the preliminary tests~(i) and~(ii) of Appendix~\ref{subsec:matching} (for $L_6$ and $L_7$, the exact tests of Appendix~\ref{subsec:edges}). All of them are excluded, except the placements of a polyhedron on itself recorded in Table~\ref{tab:same-type}.}\label{tab:coverage}
\begin{tabular}{lcrl}\toprule
source & faces of the target & placements & conclusion\\\midrule
$L_6$    &12--20&50  &$\dd(B,Q)\ge10$ at pentagons\\
$L_7$    &12--18&2   &$\dd(B,Q)\ge8$ at pentagons\\
$F_{15}$ &12--23&4206&$\dd(B,Q)\ge13$ if $B\not\cong F_{15}$\\
$F_{16a}$&12--20&690 &$\dd(B,Q)\ge10$ if $B\not\cong F_{16a}$\\
$F_{16b}$&12--20&654 &$\dd(B,Q)\ge10$ if $B\not\cong F_{16b}$\\\bottomrule
\end{tabular}
\end{table}

\begin{table}[ht]
\centering
\caption{Surviving placements of $F_{15}$, $F_{16a}$, $F_{16b}$ on themselves. Faces of the source are grouped into orbits under its symmetry group, and a dash means that no placement survives on that orbit.}\label{tab:same-type}
\begin{tabular}{lccc}\toprule
source & size of face orbit & sides & defect\\\midrule
$F_{15}$ &6&5&4\\
         &3&6&---\\
         &6&5&5\\\midrule
$F_{16a}$&4&5&6\\
         &4&6&---\\
         &4&5&6\\
         &4&5&7\\\midrule
$F_{16b}$&12&5&6\\
         &4&6&7\\\bottomrule
\end{tabular}
\end{table}

For $F_{17a}$, $F_{17b}$, $F_{17c}$ (Lemma~\ref{lem:single}(iii)) the targets have at most $18$ faces. The surviving placements of defect $7$ are all placements of the source on itself, and the faces carrying them form orbits of sizes $1,2,2,1$ for $F_{17a}$, of sizes $2,1$ for $F_{17b}$, and one orbit of size $2$ for $F_{17c}$. No placement of defect less than $7$ survives.

\subsection{Simultaneous neighbours}\label{subsec:simultaneous-data}
For $F_{15}$, each of the twelve faces in the two pentagon orbits of Table~\ref{tab:same-type} carries exactly one surviving placement of $F_{15}$ on itself, up to symmetries of the target. In each of the cyclic sequences of faces
\[
 (0,6,11,12,13,7),\qquad(3,8,9,10,4,14),
\]
consecutive faces, including the last and the first, share an edge at one of whose endpoints the two placements prescribe transverse edges of different lengths. These two $6$-cycles in $\Omega_{F_{15}}$ give Lemma~\ref{lem:simultaneous}(i).

For $F_{16a}$ there are sixteen surviving placements of $F_{16a}$ on itself, listed in Table~\ref{tab:states}. Each row gives the source face, the target face, and the images of the cyclically ordered neighbours of the source face. The four sets
\[
 \{0,5,6,7\},\qquad\{1,8,9,10\},\qquad\{2,3,4,14\},\qquad\{11,12,13,15\}
\]
are cliques in $\Omega_{F_{16a}}$, which gives Lemma~\ref{lem:simultaneous}(ii).

\begin{table}[ht]
\caption{The sixteen surviving placements of $F_{16a}$ on itself.}\label{tab:states}
\centering
\begin{tabular}{rrrl}\toprule
placement & source face & target face & images of the neighbours\\\midrule
0&0&0&$2,7,6,5,14$\\
1&1&0&$2,7,6,5,14$\\
2&3&3&$2,14,4,9,8$\\
3&3&3&$14,2,8,9,4$\\
4&4&4&$3,14,5,10,9$\\
5&6&3&$2,8,9,4,14$\\
6&6&3&$14,4,9,8,2$\\
7&7&4&$14,5,10,9,3$\\
8&9&4&$14,5,10,9,3$\\
9&10&3&$2,8,9,4,14$\\
10&10&3&$14,4,9,8,2$\\
11&12&4&$10,9,3,14,5$\\
12&13&3&$4,14,2,8,9$\\
13&13&3&$8,2,14,4,9$\\
14&14&0&$14,2,7,6,5$\\
15&15&0&$14,2,7,6,5$\\\bottomrule
\end{tabular}
\end{table}

For $F_{16b}$, Lemma~\ref{lem:simultaneous}(iii) is checked on all $24$ incident pairs consisting of a pentagon and a hexagon. In each case every surviving placement on the pentagon is incompatible with every surviving placement of $F_{16b}$ on the hexagon, because of disjoint enclosures for a transverse edge, and the face of the target across the common edge is a hexagon.

\subsection{Facets with twenty to twenty-seven faces}\label{subsec:terminal-data}
The ten graphs with $24$ vertices in the proof of Lemma~\ref{lem:stage2} are the entries $(24,i)$ of $\mathcal G_{24}$ with
\[
 i\in\{267,\ 3498,\ 4610,\ 4739,\ 4822,\ 5340,\ 6729,\ 6762,\ 6802,\ 6831\};
\]
each carries exactly one candidate set. The two graphs not excluded by the enclosures in the proof of Lemma~\ref{lem:stage1} are $(22,622)\cong\Gamma(\Tt)$ and $(22,626)\cong\Gamma(\Th)$.

In the graph $(25,23383)$ the nonagon is the face $2$, with cyclically ordered neighbours
\[
 (1,18,17,13,21,20,3,7,6).
\]
Its neighbours $1$, $13$ and $3$ are hexagons, with cyclically ordered neighbours
\[
 1\colon(2,6,10,5,14,18),\qquad 13\colon(2,17,16,12,22,21),\qquad 3\colon(2,20,24,4,8,7),
\]
and in each list the faces in positions two, three, five and six are pentagons, while the face in position four is a heptagon. These are the incidences used in the proof of Lemma~\ref{lem:stage2}.

\nocite{And70b,Atk09,Dav08,DJ91,Duf10,ERT12,RT00}
\bibliographystyle{amsalpha}
\bibliography{references}

\end{document}